\documentclass[11pt,leqno]{article}

\usepackage{amsmath,amsfonts,amscd,amssymb,theorem}

\long\def\comment#1\endcomment{}

\comment
\endcomment

\makeatletter
\begingroup
\gdef\th@dotted{\normalfont\itshape
  \def\@begintheorem##1##2{%
        \item[\hskip\labelsep \theorem@headerfont ##1\ ##2.]}%
\def\@opargbegintheorem##1##2##3{%
   \item[\hskip\labelsep \theorem@headerfont ##1\ ##2\ (##3).]}}
\endgroup
\makeatother

\theoremstyle{dotted}

\newtheorem{theorem}{Theorem}[section]
\newtheorem{lemma}[theorem]{Lemma}

\newtheorem{prop}[theorem]{Proposition}
\newtheorem{corr}[theorem]{Corollary}

\makeatletter
\begingroup
\gdef\th@upshape{\normalfont
  \def\@begintheorem##1##2{%
        \item[\hskip\labelsep \theorem@headerfont ##1\ ##2.]}%
\def\@opargbegintheorem##1##2##3{%
   \item[\hskip\labelsep \theorem@headerfont ##1\ ##2\ (##3).]}}
\endgroup
\makeatother

\theoremstyle{upshape}

\newtheorem{defn}[theorem]{Definition}
\newtheorem{remark}[theorem]{Remark}
\newtheorem{exa}[theorem]{Example}

\makeatletter
\renewcommand{\subsection}{\@startsection{subsection}{2}{0pt}{-3ex
plus -1ex minus -0.2ex}{-2mm plus -0pt minus
-2pt}{\normalfont\bfseries}} 
\renewcommand{\subsubsection}{\@startsection{subsubsection}{3}{0pt}{-3ex
plus -1ex minus -0.2ex}{-2mm plus -0pt minus
-2pt}{\normalfont\bfseries}} 
\makeatother

\makeatletter
\@addtoreset{equation}{section}
\makeatother
\renewcommand{\theequation}{\thesection.\arabic{equation}}

\newcommand{\cntrct}                
{\hspace{2pt}\raisebox{1pt}{\text{$\lrcorner$}}\hspace{2pt}}

\newcommand{\proof}[1][Proof.]{\smallskip\noindent{\em #1}}
\def\endproof{\hfill\ensuremath{\square}\par\medskip}

\def\eqref#1{\thetag{\ref{#1}}}

\let\latexref=\ref
\def\ref#1{{\normalfont{\latexref{#1}}}}

\newcommand{\wt}{\widetilde}

\newcommand{\idot}{{\:\raisebox{1pt}{\text{\circle*{1.5}}}}}
\newcommand{\hdot}{{\:\raisebox{3pt}{\text{\circle*{1.5}}}}}
\newcommand{\Z}{{\mathbb Z}}
\newcommand{\Nb}{{\mathbb N}}

\newcommand{\eps}{\varepsilon}
\renewcommand{\phi}{\varphi}

\def\dlim_#1{{\displaystyle\lim_{#1}}}

\newcommand{\Hom}{\operatorname{Hom}}

\newcommand{\id}{\operatorname{\sf id}}
\newcommand{\Id}{\operatorname{\sf Id}}
\newcommand{\gr}{\operatorname{\sf gr}}
\newcommand{\ad}{\operatorname{\sf ad}}

\newcommand{\A}{{\cal A}}
\newcommand{\D}{{\cal D}}
\newcommand{\C}{{\cal C}}
\newcommand{\E}{{\cal E}}

\newcommand{\I}{{\sf I}}
\newcommand{\T}{{\sf T}}

\newcommand{\hash}{\sharp}

\newcommand{\Cycl}{\operatorname{Cycl}}

\newcommand{\Cat}{\operatorname{Cat}}

\newcommand{\amod}{{\text{\rm -mod}}}
\newcommand{\proj}{{\text{\rm -proj}}}
\newcommand{\biproj}{{\text{\rm -pspf}}}

\newcommand{\bimod}{{\text{\rm -bimod}}}
\newcommand{\ppt}{{\sf pt}}

\newcommand{\lotimes}{\overset{\sf\scriptscriptstyle L}{\otimes}}

\newcommand{\M}{\mathcal{M}}
\newcommand{\N}{\mathcal{N}}

\newcommand{\Dd}{\operatorname{\sf D}}
\newcommand{\Nn}{\operatorname{\sf N}}

\newcommand{\sh}{\operatorname{\sf sh}}
\newcommand{\go}{{\geq 0}}
\newcommand{\B}{{\sf B}}

\newcommand{\Rr}{\rho}
\newcommand{\Ll}{\lambda}

\newcommand{\dd}{\Delta^o\Delta}

\newcommand{\g}{\mathfrak{g}}

\newcommand{\rk}{\operatorname{\sf rk}}

\title{Trace theories and localization}

\author{D. Kaledin\thanks{Partially supported by RSCF, grant number
    14-21-00053, and the Dynasty Foundation award.}}

\begin{document}

\maketitle

\tableofcontents

\section*{Introduction.}

This paper arose as an attempt to understand and generalize the
fundamental localization theorem of B. Keller proved in \cite{kel}
-- namely, the theorem claiming that whenever we have a short exact
sequence of small DG categories, we have a long exact sequence of
their Hochschild homology groups. The argument in \cite{kel} is very
well organized and its parts fit together very tightly; at first, it
is not clear how to move them or indeed whether it is possible at
all. However, on some reflection, the essential key to the argument
seems to be the following property of Hochschild homology of
associative algebras: if we have two algebras $A$, $B$, an $A
\otimes B^o$-module $M$, and a $B \otimes A^o$-module $N$, then
there is a canonical identification
\begin{equation}\renewcommand{\theequation}{*}\label{st}
HH_\idot(A,M \lotimes_B N) \cong HH_\idot(B,N \lotimes_A M).
\end{equation}
This isomorphism admits a rather straightforward generalization to
DG categories and DG bimodules, and having done this, one can
condense the crucial part of Keller's proof to a half-a-page
argument (for the convenience of the reader, I do this at the end of
the paper, in the beginning of Subsection~\ref{loc.subs}). The
isomorphism \thetag{*} itself is the derived version of the
corresponding isomorphism
$$
HH_0(A,M \otimes_B N) \cong HH_0(B,N \otimes_A M)
$$
in degree $0$, and this isomorphism is obvious from the
definitions. Thus the functor $HH_0(-,-)$, -- and consequently, its
derived functor $HH_\idot(-,-)$, -- has a ``trace-like property'':
applying the functor to a product of two bimodules gives the same
result no matter in what order we put the multiples.

This opens the door to possible generalizations, and indeed, one
such was attemped in \cite{K2}. There the trace-like property
figures very prominently, and it has even been axiomatized in the
definition of a so-called ``trace-functor'' on a monoidal category.

At first, this does not lead very far -- essentially the only
example of an additive trace functor on the category of bimodules
over an algebra $A$ known to the author is $HH_0(A,-)$, and if one
imposes some minimal compatibility conditions on trace functors for
different algebras, then one can actually prove that no other
examples exist.

However, and somewhat surprisingly, we do have non-trivial examples
of {\em non-additive} trace functors on bimodule categories -- in
fact already on the category of vector spaces over a field -- and
quite a large part of the theory can be made to work perfectly well
in the non-additive setting. This is the subject of the present
paper.

\bigskip

The paper is organized as follows. In Section~\ref{cat.sec}, we set
up the basic categorical preliminaries for the story -- we discuss
monoidal categories, algebras and modules in monoidal categories,
$2$-categories, their cyclic nerves and suchlike. This we do at
some length, and maybe in too much detail for our eventual purpose;
however, I felt that it is better to write down my take on this
story once and for all. In Section~\ref{tr.fun.sec}, we introduce
trace functors and a related notion of a ``trace theory'', and we
show how to twist Hochschild and cyclic homology by plugging a trace
functor into it. We also provide some examples; in particular, we
introduce one non-trivial non-additive example that will be the
fundamental testing case for all the later constructions (this is
Example~\ref{cycl.exa}).

Then we move on to the DG world. In Section~\ref{N.D.sec}, we recall
the basics about Dold-Kan equivalences, Dold derived functors of
non-additive functors, and so on. Since we want to consider
complexes that can have non-trivial components in both positive and
negative homological degrees, the classic theory of Dold-Puppe is
not sufficient, and we have to work with simplicial-cosimplicial
objects. At this point, in order to get a good homotopy-invariant
derived functors, we need to impose a condition on the trace
functors that we consider. We call the good functors {\em
  balanced}. We also prove some criteria for a functor to be
balanced. Fortunately, the condition is rather weak, so that
balanced functors are easy to come by (and in particular, the
functors of Example~\ref{cycl.exa} are balanced). Then in
Section~\ref{dg.alg.sec}, we define and study twisted Hochschild
homology for DG algebras and prove a derived version of the trace
isomorphisms for our twisted Hochschild homology.

In Section~\ref{dg.cat.sec}, we reap the benefits: as it happens,
the derived trace isomorphisms by themselves are enough to deduce
many of the desired applications. In particular, nothing else is
needed to extend things from DG algebras to DG categories, and to
establish derived Morita-invariance of the resulting twisted
Hochschild homology theory.

Finally, in Subsection~\ref{loc.subs}, we turn to localization. We
discover that while Keller's original proof can be restated very
concisely in our language, an analogous localization theorem for a
general trace functor is definitely not true, and one needs to
impose further conditions. At this point, we give up the attempts to
axiomatize the situation, and finish the paper by proving that at
least for the non-additive trace functor of Example~\ref{cycl.exa} --
that is, in our testing case -- the localization theorem does hold.

\bigskip

It should be emphasized that most if not all of the technical
material in this paper is not new; what I believe to be new is the
assembly that leads to the final result. In particular, I am
certainly not the first person in the world to study
categorifications of the notion of the trace -- in fact, there were
papers concerned with this already in the 1960ies, and the story
continues to attract attention. Among recent papers, one should
mention a solid and detailed theory developed in \cite{pon} -- the
notion of a ``shadow'' introduced there is very close to what I call
a trace theory.

There are also several recent papers that study the notion of a
trace in a homotopical or higher-categorical setting. I should
explain, however, that this is quite orthogonal to what I do in this
paper. Indeed, a homotopical approach would be to start with a trace
functor defined on the category of complexes of vector spaces, and
study its exactness properties by modern techniques of homotopical
algebra (model categories, infinity-categories and suchlike). The
problem is, I do not know whether trace functors like this appear
naturally. Non-additive trace functors on the category of vector
spaces do appear naturally, but since they are non-additive,
extending them to complexes is a separate and non-trivial
exercize. The only consistent way to do this known to me is the old
approach of Dold-Puppe and Kan that has nothing to do with
homotopical algebra in the modern post-Quillen sense. Of course,
having thus extended a trace functor from vector spaces to
complexes, one can apply the homotopical techniques. But in
practice, this seems pointless -- I do not see any additional
statements one might be able to prove in this way.

The trick of using simpicial-cosimplicial objects and second
quadrant bicomplexes that is the basis of Section~\ref{N.D.sec} is
also not new -- it appeared before at least in \cite{I}, and
probably elsewhere, too. However, I do not know whether the property
of a functor to be balanced has been formulated and studied
before. Neither do I know any intrinsic explanation for this
property; at this point, it looks like a rather mysterious technical
gadget.

Finally, let me emphasize again that the essential argument for the
localization theorem presented in this paper is precisely the
original argument of Keller. The only thing I did was to rearrange
its parts, to open holes in it, so that generalizations can be
plugged in.

\subsection*{Acknowledgements.} Throughout this work, I have
benefited from conversations with many people; I am particularly
grateful to Sasha Efimov, Sasha Kuznetsov and Bertrand
To\"en. Special thanks are due to Bertrand for convincing me that I
should not be ashamed of having category theory as my favourite
mathematical tool. An early version of the paper has been presented
in 2012 at a conference in Huatulco, Mexico, at a kind invitation of
Ernesto Lupercio. Substantial part of the work has been done during
my visit to University of Campinas. I am very grateful to Elizabeth
Gasparim for inviting me, and to Unicamp in general for its
wonderful working atmosphere. I am also grateful to the referee for
a very detailed and thoughtful report, with many useful corrections
and suggestions.

\section{Categorical preliminaries.}\label{cat.sec}

\subsection{Segal categories.}\label{seg.subs}

We begin by recalling a description of monodial categories and
$2$-categories in the spirit of G. Segal (this is a tiny fraction of
what is nowadays available in the literature, but we prefer to spell
things out to fix terminology and notation). 

\medskip

We need the machinery of fibered and cofibered categories of
\cite{SGA}; we summarize it as follows.

\begin{defn}
A morphism $f:c' \to c$ in a category $\C$ is {\em cartesian} with
respect to a functor $\phi:\C \to \E$ if any $f':c'' \to c$ such
that $\phi(f) = \phi(f')$ uniquely factorizes as $f' = f \circ f_0$
with $\phi(f_0)=\id$. A {\em cartesian lifting} of a morphism $f:e'
\to e$ in $\E$ is a morphism $f'$ in $\C$ cartesian with respect to
$\phi$ and such that $\phi(f')=f$. A functor $\phi:\C \to \E$ is a
    {\em prefibration} if for any $c \in \C$, any morphism $f:e' \to
    e=\phi(c)$ in $\E$ admits a cartesian lifting $f':c' \to c$. A
    prefibration is a {\em fibration} if the composition of
    cartesian morphisms is cartesian. A morphism is {\em
      cocartesian} if it is cartesian as a morphism in the opposite
    category $\C^o$ with respect to the opposite functor $\phi^o:\C^o
    \to \E^o$, and a functor $\phi$ is a {\em cofibration} if the
    opposite functor $\phi^o$ is a fibration.
\end{defn}

If $\phi:\C \to \E$ is a fibration, then cartesian liftings of
morphisms in $\E$ are defined by their targets up to a canonical
isomorphism, and sending $c \in \C$ to the source $c'$ of the
cartesian lifting $f'$ of a morphism $f:e' \to e=\phi(c)$ gives a
functor $f^*:\C_{e} \to \C_{e'}$ between the fibers $\C_{e} =
\phi^{-1}(e)$, $\C_{e'} = \phi^{-1}(e')$ of the fibration
$\phi$. The functors $f^*$ are called {\em transition functors} of
the fibration. The fibers $\C_e$, $e \in \E$ of a fibration $\phi:\C
\to \E$ together with transition functors between them patch
together to define what Grothendieck calls a contravariant {\em
  pseudofunctor} from $\E$ to the category of categories; giving a
category fibered over $\E$ is equivalent to giving such a
pseudofunctor (this is sometimes called the {\em Grothendieck
  construction}). Replacing the fibers $\C_e$ with opposite
categories $\C_e^o$ gives another pseudofunctor; the corresponding
fibered category has the same objects as $\C$, with morphisms from
$c'$ to $c$ given by isomorphism classes of diagrams
$$
\begin{CD}
c' @<{v}<< \wt{c} @>{h}>> c
\end{CD}
$$
with $\phi(v) = \id$ and cartesian $h$. Composition of morphisms is
given by taking appropriate pullbacks; these exist automatically for
any fibration $\phi:\C \to \E$. We will denote this category by
$(\C/\E)^o$.

\medskip

Now, as usual, we denote by $\Delta$ the category of finite
non-empty totally ordered sets, and we denote by $[n] \in \Delta$
the set with $n$ elements. Denote by $s,t \in [n]$ the initial
resp. the terminal element of the ordinal $[n]$. For every $[n],[m]
\in \Delta$, we have a cocartesian square in $\Delta$ of the form
\begin{equation}\label{del.sq}
\begin{CD}
[1] @>{s}>> [n]\\
@V{t}VV @VVV\\
[m] @>>> [n+m-1].
\end{CD}
\end{equation}

\begin{defn}\label{spec.defn}
A fibration $\pi:\C^\flat \to \Delta$ is {\em special} if
\begin{enumerate}
\item the fiber $\C^\flat_{[1]}$ over $[1] \in \Delta$ is discrete,
  and
\item for any square \eqref{del.sq}, the corresponding diagram
$$
\begin{CD}
\C^\flat_{[1]} @<{s^*}<< \C^\flat_{[n]}\\
@A{t^*}AA @AAA\\
\C^\flat_{[m]} @<<< \C^\flat_{[n+m-1]}
\end{CD}
$$
of fibers of the fibration $\C^\flat \to \Delta$ and transition
functors between them is weakly cartesian -- that is, it induces an
equivalence between $\C^\flat_{[n+m-1]}$ and the category of pairs
$\langle a,b \rangle$ of objects $a \in \C^\flat_{[n]}$, $b \in
\C^\flat_{[m]}$ equipped with an isomorphism $s^*a \cong t^*b$.
\end{enumerate}
A {\em $2$-category $\wt{\C}$} is the pair of a category
$\Nn(\wt{\C})$ called the {\em nerve} of $\C$ and a special
fibration $\Nn(\wt{\C}) \to \Delta$.
\end{defn}

The usual definition of a $2$-category can be recovered as follows.
The objects are objects of $\Nn(\wt{\C})_{[1]}$. For any two objects
$c,c' \in \Nn(\wt{\C})_{[1]}$, the category $\C(c,c')$ of morphisms
from $c$ to $c'$ is the category of diagrams $c \to \wt{c} \gets c'$
in $\Nn(\wt{C})$ such that the functor $\pi$ sends it to the diagram
$$
\begin{CD}
[1] @>{s}>> [2] @<{t}<< [1]
\end{CD}
$$
in $\Delta$. It is well-known that this establishes a one-to-one
correspondence between $2$-categories in the usual sense and in the
sense of Definition~\ref{spec.defn}.

\begin{exa}
For any small category $\Phi$, its nerve is a simplicial set, and it
corresponds by the Grothendieck construction to a discrete fibration
over $\Delta$; we denote it by $\Nn(\Phi)/\Delta$. This fibration is
special, and it corresponds to $\Phi$ considered as a $2$-category.
\end{exa}

\begin{exa}
Unital associative monoidal categories $\C$ are in one-to-one
correspondence with $2$-categories $B(\C)$ with one object $o$; $\C$
is the category of endofunctors $B(\C)(o,o)$. In the context of
Definition~\ref{spec.defn}, these correspond to special fibrations
$\C^\flat/\Delta$, $\C^\flat = \Nn(B(\C))$ such that $\C^\flat_{[1]}
= \ppt$ is the the point category.
\end{exa}

\begin{defn}\label{2fun.defn}
A {\em functor} $\phi:\wt{\C}_1 \to \wt{\C}_2$ between
$2$-categories is a cartesian functor $\Nn(\phi):\Nn(\wt{\C}_1) \to
\Nn(\wt{\C}_2)$ -- that is, a functor that commutes with projections
to $\Delta$ and sends cartesian maps to cartesian maps. A {\em
  morphism} between functors $\phi$, $\phi'$ is a morphism
$\Nn(\phi) \to \Nn(\phi')$.
\end{defn}

\begin{remark}
Of course one would expect $2$-categories to form a $3$-category, so
that we have morphisms between morphisms between functors. But we
will not need this.
\end{remark}

\begin{exa}
For any two unital associative monoidal categories $\C_1$, $\C_2$, a
functor $B(\C_1) \to \B(\C_2)$ is the same things as a unital
monoidal functor $\C_1 \to \C_2$.
\end{exa}

For any $2$-category $\wt{\C}$, one can define a ``half-opposite''
$2$-category $\wt{\C}^\natural$ by keeping the same objects, and
replacing the morphism categories with the opposite categories. On
the level of nerves, this corresponds to replacing a fibration
$\C^\flat/\Delta$ with the fibration
$(\C^\flat/\Delta)^o/\Delta$. Repeating this procedure twice, we
recover our original fibration $\C^\flat/\Delta$. In particular,
giving a unital associative monoidal structure on a category $\C$ is
completely equivalent to giving such a structure on the opposite
category $\C^o$. It will be often more convenient for us to consider
the $2$-category $B(\C^o) = B(\C)^\natural$ instead of $B(\C)$. We
will call its nerve $\C^{o\flat}/\Delta$ the {\em Segal fibration}
associated to $\C$.

\begin{defn}
A map $f:[n'] \to [n]$ in $\Delta$ is an {\em anchor map} if it
identifies $[n']$ with the subset $[a,b] \subset [n]$ of elements $v
\in [n]$ such that $a \leq v \leq b$, for some $a,b \in [n]$. For
any fibration $\phi:\C \to \Delta$, a map $f$ in $\C$ is an anchor
map if it is cartesian with respect to $\phi$ and $\phi(f)$ is an
anchor map in $\Delta$.
\end{defn}

\begin{defn}\label{ps.fun.def}
A {\em pseudofunctor} $\phi:\wt{\C}_1 \to \wt{\C}_2$ between
$2$-categories is a functor $\Nn(\phi):\Nn(\wt{\C}_1) \to
\Nn(\wt{\C}_2)$ that commutes with projections to $\Delta$ and sends
anchor maps to anchor maps.
\end{defn}

We note that this is weaker than the notion of a pseudofunctor used
in \cite{SGA}. In particular, for two monoidal categories $\C$,
$\C'$, a pseudofunctor $\phi:B(\C) \to B(\C')$ corresponds to a {\em
  pseudotensor} functor $\overline{\phi}:\C \to \C'$ -- that is, a
functor $\overline{\phi}$ equipped with a functorial map
$$
\overline{\phi}( - \otimes - ) \to \overline{\phi}(-) \otimes
\overline{\phi}(-),
$$
subject to obvious associativity and unitality constraints. The map
need not be an isomorphism.

\subsection{Enrichments.}\label{enr.subs}

Assume given a unital associative monoidal category $\C$ with tensor
product $- \otimes -$, and let $\C^{o\flat}/\Delta$ be the
corresponding Segal fibration.

\begin{defn}\label{enr.def}
\begin{enumerate}
\item A {\em $\C$-enrichment} for a small category $\Phi$ is a
  pseudofunctor $\phi:\Nn(\Phi) \to \C^{o\flat}$. An
  enrichment is called {\em precise} if it is cartesian.
\item For any two small categories $\Phi_1$, $\Phi_2$ with
  $\C$-enrichments $\phi_1$, $\phi_2$, a {\em $\C$-enrichment for a
  functor} $F:\Phi_1 \to \Phi_2$ is a morphism $\phi_2 \circ \Nn(F)
  \to \phi_1$. A $\C$-enrichment for $F$ is {\em precise} if it is
  an isomorphism $\phi_2 \circ \Nn(F) \cong \phi_1$.
\end{enumerate}
Composition of enriched functors is defined in an obvious way.
\end{defn}

Explicitly, an enrichment is given by associating an object $\phi(f)
\in \C$ to any morphism $f$ in $\Phi$, and a composition map
\begin{equation}\label{compopo}
\phi(f) \otimes \phi(g) \to \phi(f \circ g)
\end{equation}
to any composable pair of morphisms $f$, $g$ (were we to use
$\C^\flat$ instead of $\C^{o\flat}$, the map would go in the
opposite direction). These data are related by all sorts of higher
associativity and unitality isomorphisms packed into the functor
$\phi$ of Definition~\ref{enr.def}. An enrichment is precise if all
the maps \eqref{compopo} are isomorphisms (this is rare in practice,
but we will need it). Analogously, an enrichment for a functor $F$
is given by maps
$$
\phi_1(f) \to \phi_2(F(f)),
$$
one for every morphism $f$ in $\Phi_1$, subject to various
constraints. It is precise if all maps are isomorphisms.

\begin{remark}\label{usual.enr}
Our notion of enrichment is slightly more general then the usual
one; to obtain the usual notion of a $\C$-enriched small category
with the set of object $S$, let $E(S)$ be the category with $S$ as
the set of objects and exactly one morphism between every two
objects, and apply Definition~\ref{enr.def} to $\Phi = E(S)$.
\end{remark}

\begin{exa}\label{alg.exa}
A $\C$-enrichment of the point category $\ppt$ is the same thing as
a unital associative algebra object in $\C$.
\end{exa}

The nerve $\Nn(\ppt)$ of the point category is the category $\Delta$
itself, so that an enrichment for $\ppt$ is the same thing as a
section $\Delta \to \C^{o\flat}$ of the projection $\C^{o\flat} \to
\Delta$ that is a pseudofunctor -- that is, sends anchor maps to
anchor maps.

\begin{exa}\label{mon.exa}
More generally, assume given a monoid $M$, and let $B(M)$ be the
category with one object and endomorphism monoid $M$. Then a
$\C$-enrichment for $B(M)$ is the same thing as an $M$-graded
algebra in $\C$ -- that is, a collection of objects $A_m \in \C$, $m
\in M$, and multiplication maps $A_m \otimes A_{m'} \to A_{m \cdot
  m'}$, subject to obvious associativity and unitality conditions.
\end{exa}

\begin{exa}\label{mod.exa}
Let us associate a small category to any partially ordered set in
the usual way, and let $[2]_\Delta$ be the category corresponding to
the totally ordered set $[2] \in \Delta$. Then a $\C$-enrichment for
$[2]_\Delta$ is the same thing as an algebra object $A_1$ in $\C$,
an algebra object $A_2$ in $\C$, and a left module $M$ over $A_1^o
\otimes A_2$, where $A^o$ is the opposite algebra.
\end{exa}

Unfortunately, one cannot define $A$-bimodules in a similar
fashion. It is clear from Example~\ref{mod.exa} that pairs $\langle
A,M \rangle$ naturally correspond to functors
$$
\Sigma/\Delta \to \C^{o\flat}
$$
sending anchor maps to anchor maps, where the fibered category
$\Sigma/\Delta$ corresponds by the Grothendieck construction to the
simplicial circle -- that is, to the coequalizer of the diagram
$$
[1] \begin{array}{c}\overset{s}{\longrightarrow}\\[-2mm]
  \underset{t}{\longrightarrow}\end{array}
[2]
$$
in the category of simplicial sets. However, $\Sigma/\Delta$ is not
a special fibration in the sense of Definition~\ref{spec.defn}. Thus
to work with bimodules, we use a workaround. First, we assume that
$\C$ has a terminal object $0$ such that $0 \otimes M = 0$ for any
$M \in \C$ (this is harmless, since such a terminal object can
always be formally added to $\C$). Next, we observe that the
category of pairs $\langle A,M \rangle$ of an algebra $A$ and an
$A$-bimodule $M$ is equivalent to the category algebras $A_\idot$ in
$\C$ graded by the monoid $\Nb$ of non-negative integers such that
$A_i = 0$ for $i \geq 2$. The equivalence sends $A_\idot$ to
$\langle A_0,A_1 \rangle$. On the simplicial level, we have a
natural embedding $\Sigma \to \Nn(B(\Nb))$, where $B(\Nb)$ is as in
Example~\ref{mon.exa}, every $\C$-enrichment for $B(\Nb)$ induces a
functor $\Sigma \to \C^{o\flat}$ by restriction, and we use the
right-adjoint functor to this restriction to promote bimodules to
$\Nb$-graded algebra objects.

\subsection{Cyclic categories.}

To proceed further, we need to recall the basics about A. Connes'
{\em cyclic category $\Lambda$}. Let $[\infty]$ be $\Z$ considered
as a partially ordered set, with its usual order. Let
$\sigma:[\infty] \to [\infty]$ be the map sending $a \in \Z$ to
$a+1$. By definition, $\Lambda_\infty$ is the category whose objects
$[n]$ are numbered by positive integers $n \geq 1$, and whose
morphisms $f \in \Lambda_\infty([n],[m])$ are maps $f:[\infty] \to
[\infty]$ such that $f \circ \sigma^n = \sigma^m \circ f$. For
any such $f$, $\sigma(f)$ given by $\sigma(f)(a) = \sigma^m(f(a)) =
f(\sigma^n(a))$ also defines a morphism in $\Lambda_\infty$, and for
any two composable morphisms $f$, $g$, we have $\sigma(f) \circ g =
f \circ \sigma(g) = \sigma(f \circ g)$. The category $\Lambda$ has
the same objects $[n]$ as $\Lambda_\infty$, and the morphisms given
by
$$
\Lambda([n],[m]) = \Lambda_\infty([n],[m])/\sigma.
$$
Note that the category $\Lambda_\infty$ is self-dual: the
equivalence $\Lambda_\infty \cong \Lambda_\infty^o$ is identical on
objects, and sends $f \in \Lambda_\infty([n],[m])$ to $f^o \in
\Lambda_\infty([m],[n])$ given by
$$
f^o(a) = \min\{b \in \Z \mid f(b) \geq a\}.
$$
This descends to an equivalence $\Lambda \cong \Lambda^o$.

For any $[n] \in \Lambda$, let $V([n]) = \Z/n\Z$; then by
definition, $V$ is functor from $\Lambda$ to the category of finite
sets (in fact, $V([n]) = \Lambda([1],[n])$, so that the functor $V$
is representable). One can extend $V$ to a functor from $\Lambda$ to
the category $\Cat$ of small categories by considering the small
category $[\infty]_\Delta$ corresponding to the partially ordered
set $[\infty]$, and letting $[n]_\Lambda$ to be the quotient
category $[\infty]_\Delta/\sigma^n$ for any $n \geq 1$. This means
that objects of $[n]_\Lambda$ are elements $v \in V([n]) = \Z/n\Z$,
and morphisms are given by
$$
\Hom_{[n]_\Lambda}(v,v') = \coprod_{l \in
  \Z}\Hom_{[\infty]_\Delta}(a,a' + nl)
$$
for any $v,v' \in \Z/n\Z$ represented by some $a,a' \in
\Z$. Equivalently, morphisms from $v$ to $v'$ are given by
non-negative integers $l \geq 0$ such that $l = v'-v \mod n$. We
denote the morphism $v \to v+l$ corresponding to $l \geq 0$ by
$l_v$, and for any morphism $f = l_v$ in $[n]_\Lambda$, we denote $l
= l(f)$ and call it the {\em length} of $f$ (note that morphisms of
length $0$ are identity morphisms, $0_v = \id_v$). Sending $[n]$ to
$[n]_\Lambda$ gives a faithful embedding $\Lambda \subset \Cat$.

\begin{remark}\label{path}
One can visualize an object $[n] \in \Lambda$ as a wheel quiver with
$n$ vertices $v \in V([n])$ and $n$ edges oriented clockwise. The
category $[n]_\Lambda$ is the path category of this quiver, with
length of a morphism being the length of the corresponding path.
Morphisms in $\Lambda$ do not all correspond to quiver maps, but
they do correspond to functors between the path categories.
\end{remark}

For any map $f:[n'] \to [n]$ in $\Lambda$ represented by a
order-preserving map $\wt{f}:[\infty] \to [\infty]$, and any $v \in
V([n])$, the preimage $f^{-1}(v) \subset V([n'])$ is naturally
identified with the preimage $\wt{f}^{-1}(\wt{v}) \subset [\infty]$
of any representative $\wt{v} \in [\infty]$ of the element
$v$. Therefore it inherits a natural total order. In particular, let
$\Lambda/[1]$ be the category of objects $[n] \in \Lambda$ equipped
with a map $f:[n] \to [1]$. Then sending to the preimage of the only
element in $V([1])$ gives a natural functor $\Lambda/[1] \to
\Delta$. One easily checks that this functor is an equivalence of
categories. Applying duality, we obtain an equivalence
$[1]\backslash\Lambda \cong (\Lambda/[1])^o \cong \Delta^o$, where
$[1]\backslash\Lambda$ is the category of objects $[n] \in \Lambda$
equipped with a map $[1] \to [n]$. Composing it with the forgetful
functor $[1]\backslash\Lambda \to \Lambda$, we obtain a natural
functor
\begin{equation}\label{j}
j:\Delta^o \to \Lambda.
\end{equation}
On objects, $j$ sends $[n] \in \Delta$ to $[n] \in \Lambda$ (this
explains our non-standard enumeration of simplices).

More generally, let $\Lambda/[n]$ be the category of objects $[n']
\in \Lambda$ equipped with a map $f:[n'] \to [n]$, and let
$\Lambda/_s[n] \subset \Lambda/[n]$ be the full subcategory spanned
by surjective $f:[n'] \to [n]$. Then sending such an $f$ to the
collection $f^{-1}(v)$, $v \in V([n])$ gives a functor
$$
\Lambda/_s[n] \to \Delta^{V([n])} = \Delta^n.
$$
This functor is also an equivalence of categories; composing the
inverse equivalence with the natural embedding gives a functor
\begin{equation}\label{j.n}
j_n:\Delta^n \to \Lambda/[n].
\end{equation}

\subsection{Cyclic nerves.}

Assume now given a $2$-category $\wt{\C}$ in the sense of
Definition~\ref{spec.defn}.

\begin{defn}\label{cycl.nerv.def}
The {\em cyclic nerve} $\Lambda(\wt{\C})$ of the $2$-category
$\wt{\C}$ is the category of pairs $\langle [n],\phi \rangle$ of an
object $[n] \in \Lambda$ and a functor $\phi:[n]_\Lambda \to
\wt{\C}$ in the sense of Definition~\ref{2fun.defn}, with morphisms
given by pairs of a morphism $f:[n'] \to [n]$ and an isomorphism $f
\circ \phi' \cong \phi$.
\end{defn}

Sending $\langle [n],\phi \rangle$ to $[n]$ gives a functor
\begin{equation}\label{cycl.fib}
\Lambda(\wt{\C}) \to \Lambda.
\end{equation}
This functor is obviously a fibration. If $\wt{\C} = \Phi$ is a
small $1$-category, so that its nerve $\Nn(\Phi)/\Delta$ is discrete,
then the fibration \eqref{cycl.fib} is discrete, since $\Nn(\Phi)$
has no non-trivial isomorphisms. Thus the cyclic nerve
$\Lambda(\Phi)$ is effectively a functor from $\Lambda^o$ to
sets. Its value on an object $[n] \in \Lambda$ is the set of
diagrams $\phi_1 \to \dots \to \phi_n \to \phi_1$ in the category
$\Phi$.

Another particular case that we will need is $\wt{\C} = B(\C)$ for a
unital associative monoidal category $\C$. In this case, we denote
\begin{equation}\label{c.hash}
\C^\hash = \Lambda(B(\C^o))^o,
\end{equation}
so that \eqref{cycl.fib} induces a cofibration $\C^\hash \to
\Lambda^o \cong \Lambda$. Objects of $\C^\hash$ correspond to pairs
$\langle [n],\phi \rangle$ of an object $[n] \in \Lambda$ and a
precise $\C^o$-enrichment of the category $[n]_\Lambda$ in the sense
of Definition~\ref{enr.def}. As in Remark~\ref{path}, the category
$[n]_\Lambda$ is the path category of the wheel quiver with $n$
vertices and $n$ edges, so that to giving such a precise
$\C^o$-enrichment is equivalent to giving a collection $M_\idot$ of
$n$ objects $M_0,\dots,M_{n-1} \in \C^o$ corresponding to maps $v
\to v+1$, $v \in \Z/n\Z$ represented by the integer $1$. Therefore
the fiber $\C^\hash_{[n]}$ of the cofibration $\C^\hash \to \Lambda$
is naturally identified with the product $\C^{V([n])}$ of copies of
$\C^o$ numbered by elements $v \in V([n])$. The transition functor
$f_*:\C^{V([n'])} \to \C^{V([n])}$ of the cofibration
$\C^\hash/\Lambda$ corresponding to a map $f:[n'] \to [n]$ is the
given by
\begin{equation}\label{C.hash.eq}
(f_*(M_\idot))_v = \bigotimes_{v' \in f^{-1}(v)}M_{v'}
\end{equation}
for any $v \in V([n])$ and $M_\idot = \{M_{v'}|v' \in V([n'])\} \in
\C^{V([n])}$ (where the tensor product is taken in the natural order
on the set $f^{-1}(v)$, and if the set $f^{-1}(v)$ is empty, the
component $(f_*(M_\idot))_v$ of $f_*(M_\idot) \in \C^{V([n])}$ is
the unit object $1 \in \C$).

\medskip

The construction of the cyclic nerve is obviously functorial with
respect to functors between $2$-categories -- any functor $\wt{\C}_1
\to \wt{\C}_2$ induces a functor $\Lambda(\wt{\C}_1) \to
\Lambda(\wt{\C}_2)$ cartesian with respect to \eqref{cycl.fib}. We
finish this section by observing that $\Lambda(-)$ is in fact even
functorial with respect to pseudofunctors in the sense of
Definition~\ref{ps.fun.def}.

\medskip

To see this, for any integers $a,b \in [\infty]$, $a \leq b$, let
$[a,b] \subset [\infty]$ be the subset of all $v \in [\infty]$ with
$a \leq v \leq b$. Let $\I$ be the partially ordered set of all
these interval subsets $[a,b] \subset [\infty]$ ordered by inclusion
(equivalently, $\I$ is the sets of pairs $\langle a,b \rangle$, $a,b
\in [\infty]$, $a \leq b$, with the order given by $\langle a,b
\rangle \leq \langle a',b' \rangle$ iff $a' \leq a \leq b \leq
b'$). Then all the intervals $[a,b]$ are finite non-empty totally
ordered sets, so that we have a natural functor $e:\I \to
\Delta$. We also have
$$
[\infty]_\Delta = \lim_{\overset{\I}{\to}}e(i)_\Delta,
$$
so for that for any fibration $\C^\flat/\Delta$, giving a cartesian
functor $\Nn([\infty]_\Delta) \to \C^\flat$ is equivalent to giving
a compatible system of cartesian functors $\Nn(e(i)_\Delta) \to
\C^\flat$, $i \in \I$. Since for any $i \in \I$, the nerve
$\Nn(e(i)_\Delta)$ is the functor represented by $e(i) \in \Delta$,
this is in turn equivalent to giving a cartesian section
\begin{equation}\label{c.eq}
c:\I \to e^*\C^\flat
\end{equation}
of the induced fibration $e^*\C^\flat/\I$.

However, by the very definition of the functor $e$, it sends every
map in $\I$ to an anchor map in $\Delta$. Therefore if we have two
fibrations $\C_1^\flat/\Delta$, $\C^\flat_2/\Delta$, and a functor
$\phi:\C_1^\flat \to \C^\flat_2$ that commutes with projections to
$\Delta$ and sends anchor maps to anchor maps, then for any
cartesian section $s:\I \to e^*\C_1^\flat$ of \eqref{c.eq}, the
section $\phi \circ s:\I \to e^*\C_2^\flat$ is also cartesian.

Now, for every $[n] \in \Lambda$, we have $\Nn([n]_\Lambda) =
\Nn([\infty]_\Delta)/\sigma^n$, so that a cartesian functor
$\Nn([n]_\Lambda) \to \C^\flat$ is the same thing as a
$\sigma^n$-equivariant cartesian functor $\Nn([\infty]_\Delta) \to
\C^\flat$. These therefore correspond to $\sigma^n$-equivariant
cartesian sections \eqref{c.eq}. Then for any $2$-categories
$\wt{\C}_1$, $\wt{\C_2}$ and a pseudofunctor $\phi:\wt{\C}_1 \to
\wt{\C}_2$, composing $\sigma^n$-equivariant cartesian sections
\eqref{c.eq} with $\phi$ defines a canonical functor
\begin{equation}\label{la.phi}
\Lambda(\wt{\C}_1) \to \Lambda(\wt{\C}_2).
\end{equation}
This functor commutes with projections \eqref{cycl.fib} to
$\Lambda$, although it is not neccesarily cartesian with respect to
these projections.

\section{Trace functors.}\label{tr.fun.sec}

\subsection{Definitions.}\label{tr.def.subs}

Assume given an associative monoidal category $\C$ with a unit
object $1 \in \C$, and a category $\E$. Denote by $\otimes$ the
product in $\C$.

\begin{defn}\label{tr.fu.def}
A {\em trace functor} from $\C$ to $\E$ is a collection of a functor
$F:\C \to \E$ and isomorphisms
\begin{equation}\label{M.N}
\tau_{M,N}:F(M \otimes N) \cong F(N \otimes M)
\end{equation}
for any $M,N \in \C$ such that
\begin{enumerate}
\item $\tau_{M,N}$ are functorial in $M$ and $N$,
\item for any object $M \in \C$, we have $\tau_{1,M} = \id$, and
\item for any three objects $M,N,L \in \C$, we have
\begin{equation}\label{M.N.L}
\tau_{L,M,N} \circ \tau_{N,L,M} \circ \tau_{M,N,L} = \id,
\end{equation}
\end{enumerate}
where in \thetag{ii}, we use the unitality isomorphisms to identify
$1 \otimes M \cong M \cong M \otimes 1$, and in \thetag{iii},
$\tau_{A,B,C}$ for any $A,B,C \in \C$ is the composition of the map
$\tau_{A,B \otimes C}$ and the map induced by the associativity
isomorphism $(B \otimes C) \otimes A \cong B \otimes (C \otimes A)$.
\end{defn}

\begin{exa}\label{triv.exa}
If $\C$ is a symmetric monoidal category, then every functor $F:\C
\to \E$ has a tautological structure of a trace functor, with
$\tau_{M,N}$ induced by the commutativity morphisms in $\C$.
\end{exa}

\begin{exa}
If $\C$ is the category of bimodules over an algebra $A$ over a
commutative ring $k$, then $0$-th Hochschild homology
$$
HH_0(A,M) = M/\{am-ma \mid a \in A,m \in M\}
$$
has a natural structure of a trace functor.
\end{exa}

\begin{exa}\label{cycl.exa}
For a non-trivial example of a trace functor on a symmetric monoidal
category, let $\C$ be the category of projective modules over a
commutative ring $k$, and fix an integer $l \geq 2$. For any $M \in
\C$, let $\sigma_M:M^{\otimes l} \to M^{\otimes l}$ be the order-$l$
permutation. For any $M,N \in \C$, let
$$
\wt{\tau}_{M,N}:(M \otimes N)^{\otimes l} \to (N \otimes M)^{\otimes
  l}
$$
be the product of the $l$-th tensor power of the commutativity
morphism and the automorphism $\sigma_N \otimes \id_{M^{\otimes l}}$
of $(N \otimes M)^{\otimes l} \cong N^{\otimes l} \otimes M^{\otimes
  l}$. Then $\wt{\tau}_{M,N}$ obviously commutes with $\sigma_{N
  \otimes M}$, thus induces a map
$$
\tau_{M,N}:(M \otimes N)^{\otimes l}_\sigma \to (N \otimes
M)^{\otimes l}_\sigma
$$
on the modules of coinvariants with respect to $\sigma$. Setting
$\Cycl^l(M) = M^{\otimes l}_\sigma$ with these structure maps
$\tau_{M,N}$ gives a trace functor from $\C$ to itself.
\end{exa}

To rephrase Definition~\ref{tr.fu.def} in the language of the Segal
machine of Subsection~\ref{seg.subs}, consider the cyclic nerve
$\C^\hash/\Lambda$ of the monoidal category $\C$, as in
\eqref{c.hash}.

\begin{lemma}\label{tr.hash}
For any category $\E$, giving a trace functor $\langle
F,\tau_{\idot,\idot} \rangle$ from $\C$ to $\E$ is equivalent to
giving a functor $F_\hash:\C^\hash \to \E$ which inverts all
cocartesian maps.
\end{lemma}

Here we say that a functor $F$ {\em inverts a map $f$} if $F(f)$ is
an invertible map.

\proof{} For any set $S$, let $E(S)$ be the groupoid whose objects
are elements of $S$, and with exactly one morphism between every two
objects. Then by definition, a functor $F:E(S) \to \E$ to some
category $\E$ is the same thing as a functor $F_0:S \to \E$, where
$S$ is considered as a discrete category, and an isomorphism
$\tau:\pi_1^*F \to \pi_2^F$, where $\pi_1,\pi_2:S \times S \to S$
are the natural projections; the isomorphism $\tau$ should restrict
to the identity isomorphism on the diagonal $S \subset S \times S$
and satisfy the obvious associativity condition on the triple
product $S \times S \times S$.

This construction is functorial in $S$. Thus if we have a discrete
cofibration $\pi:I_1 \to I$ with small fibers, it can be applied
pointwise to the fibers $\pi^{-1}(i) \subset I_1$, $i \in I$. We
obtain a category $E(\pi)$ cofibered over $I$, with fibers
identified with $E(\pi^{-1}(i))$, $i \in I$. Moreover, for any $l
\geq 2$, let
$$
I_l = I_1 \times_I \dots \times_I I_1
$$
be the fibered product of $l$ copies of $I_1$, taken in $\Cat$, and
let $\pi_1,\pi_2:I_2 \to I_1$ be the natural projections. Then
giving a functor $F:E(\pi) \to \E$ is equivalent to giving a functor
$F^1:I_1 \to \E$ equipped with an isomorphism
\begin{equation}\label{ta}
\tau:F^1 \circ \pi_1 \cong F^1 \circ \pi_2
\end{equation}
that restricts to $\Id$ on the diagonal $I_1 \subset I_2$ and
satisfies the obvious associativity condition after lifting to
$I_3$. Since the projection $E(\pi) \to I$ is an equivalence of
categories, this also describes functors $F:I \to \E$.

Let now $\Lambda^1 = \Lambda/[1]$, the category of pairs $\langle
[n],f \rangle$ of an object $[n] \in \Lambda$ and a map $f:[n] \to
[1]$, with the natural forgetful functor $\Lambda^1 \to \Lambda$
sending $\langle [n],f \rangle$ to $[n]$ (the reader might recall
that we have an the equivalence $\Lambda/[1] \cong \Delta$, but we
will not need it in this proof). Let $I = \C^\hash$, $I_1 =
\C^\hash_1 \to \Lambda^1 \times_\Lambda \C^\hash$, and let $\pi:I_1
\to I$ be the natural projection. Since the forgetful functor
$\Lambda^1 \to \Lambda$ is a discrete cofibration with small fibers,
$\pi$ is also a discrete cofibration with small fibers. Therefore
the above discussion applies -- giving a functor $F_\hash:\C_1^\hash
\to \E$ is equivalent to giving a functor $F_\hash^1:\C_1^\hash \to
\E$ and a map $\tau$ of \eqref{ta}, subject to the unitality and
associativity condition.

It remains to notice that such a functor $F_\hash$ inverts
cocartesian maps if and only if so does $F^1_\hash$; for
$F^1_\hash$, this happens if and only it factors through the natural
projection
$$
\C_1^\hash \to \C
$$
sending $M_\idot \times \langle [n],f \rangle \in \C_1^\hash$ to
$f_*(M_\idot)$. Thus $F^1_\hash$ is completely defined by a functor
$F:\C \to \E$. Then specifying the isomorphism $\tau$ of \eqref{ta}
is equivalent to specifying the isomorphisms $\tau_{M,N}$ of
\eqref{M.N}, the unitality condition is equivalent to
Definition~\ref{tr.fu.def}~\thetag{ii}, and the associativity
condition is equivalent to \eqref{M.N.L}.
\endproof

\subsection{Twisted cyclic homology.}\label{tw.ho.subs}

Assume given a unital associative algebra object $A$ in a unital
associative monoidal category $\C$, and let $\alpha:\Delta =
\Nn(\ppt) \to \C^{o\flat}$ be the corresponding $\C$-enrichment of
the point category $\ppt$, as in Example~\ref{alg.exa}. Then
$\alpha$ induces a section
$$
\Lambda(\alpha)^o:\Lambda \to \C^\hash,
$$
as in \eqref{la.phi}. We will denote $A_\hash = \Lambda(\alpha)^o$.

Explicitly, the component $A_\hash([n])_v$ of $A_\hash([n]) \in
\C^{V([n])} \cong \C^\hash_{[n]}$ is $A$ for any $[n] \in \Lambda$
and $v \in V([n])$, and for any $f:[n'] \to [n]$, the corresponding
map $A_\hash(f):A_\hash([n']) \to A_\hash([n])$ is given by
\begin{equation}\label{A.hash.f}
A_\hash(f) = \bigotimes_{v \in V([n])}m_{f^{-1}(v)},
\end{equation}
where $m_{f^{-1}(v)}:A^{\otimes f^{-1}(v)} \to A$ is the product map
in the algebra object $A$.

Moreover, assume given a bimodule $M$ over the algebra $A$, add the
terminal object $0$ to $\C$ if necessary, and let $\mu:\Nn(\Nb) \to
\C^{o\flat}$ be the corresponding $\C$-enrichment of the category
$B(\Nb)$, as in Subsection~\ref{enr.subs}. Note that we actually have
$B(\Nb) \cong [1]_\Lambda$. Then again by \eqref{la.phi}, $\mu$
induces a functor
$$
\Lambda(\mu):\Lambda([1]_\Lambda)^o \to \C^\hash
$$
commuting with the projections to $\Lambda \cong \Lambda^o$. Sending
$f:[n] \to [1]$ to $f_\Lambda:[n]_\Lambda \to [1]_\Lambda$ gives a
natural functor
\begin{equation}\label{e.eq}
\Lambda/[1] \to \Lambda([1]_\Lambda),
\end{equation}
and composing the opposite functor with $\Lambda(\mu)$ gives a
functor
\begin{equation}\label{wt.M.A.eq}
(M/A)_\hash:\Delta^o \cong (\Lambda/[1])^o \to \C^\hash,
\end{equation}
again commuting with projections to $\Lambda$ (where the projection
on the left-hand side is the functor $j:\Delta^o \to \Lambda$ of
\eqref{j}). For the diagonal bimodule $M = A$, we have
$$
(A/A)_\hash = j^*A_\hash.
$$
Explicitly, for any $[n] \in \Delta$, $V(j([n]))$ is equipped with a
distinguished element $o$; then the component $(M/A)_\hash([n])_v$
of the object $(M/A)_\hash([n]) \in \C^{V(j([n]))}$ is $M$ if $v=o$
and $A$ otherwise. For any $f:[n'] \to [n]$, $(M/A)_\hash(f)$ is
given by \eqref{A.hash.f}, with the terms $m_{f^{-1}(v)}$ induced
either by multiplication in $A$, if $v \neq o$, or by the $A$-action
on $M$ if $v=o$ (since $j(f):j([n']) \to j([n])$ by definition sends
$o$ to $o$, this is well-defined).

\medskip

Assume now given a trace functor $\langle F,\tau_{\idot,\idot}
\rangle$ from $\C$ to some category $\E$. Extend it to a functor
$F_\hash:\C^\hash \to \E$ as in Lemma~\ref{tr.hash}.

\begin{defn}\label{tw.hash.def}
The {\em $F$-twisted functor} $FA_\hash:\Lambda \to \E$ is given by
$$
FA_\hash = F_\hash \circ A_\hash,
$$
and the {\em $F$-twisted functor} $F(M/A)_\hash:\Delta^o \to \E$ is
  given by
$$
F(M/A)_\hash = F_\hash \circ (M/A)_\hash.
$$
\end{defn}

If $\C$ is the category of flat modules over a commutative ring $k$,
so that $A$ is a flat $k$-algebra, $\E$ is the category of all
$k$-modules, and $F:\C \to \E$ is the natural embedding, then the
functor $A_\hash = FA_\hash$ gives rise to the Hochschild homology
$HH_\idot(A)$ and the cyclic homology $HC_\idot(A)$ of the algebra
$A$ -- we have
$$
HC_\idot(A) = H_\idot(\Lambda,A_\hash), \qquad HH_\idot(A) =
H_\idot(\Delta^o,j^*A_\hash),
$$
where $H_\idot(\Lambda,-)$, resp. $H_\idot(\Delta^o,-)$ stands for
the homology groups of the corresponding small categories (that is,
for the derived functors of the direct limit functors
$\lim_{\overset{\Lambda}{\to}}$
  resp. $\lim_{\overset{\Delta^o}{\to}}$). These homology groups can
  be computed by well-known standard complexes (see
  e.g. \cite[Appendix]{ft}). In particular, the homology of the
  category $\Delta^o$ with coefficients in some simplicial object $E
  \in \Delta^ok\proj$ can be computed by the standard complex
  $E_\idot$ with terms $E_n = E([n+1])$ and the usual differential
  given by the alternating sum; this produces the Hochschild complex
  of the algebra $A$.  For an $A$-bimodule $M$, the simplicial
  $k$-module $(M/A)_\hash = F(M/A)_\hash$ gives Hochschild homology
  of the algebra $A$ with coefficients in $M$ -- we have
$$
HH_\idot(A,M) = H_\idot(\Delta^o,(M/A)_\hash)
$$
for any bimodule $M$ over the algebra $A$.

In the general case, if the target category $\E$ is abelian,
Definition~\ref{tw.hash.def} gives rise to twisted versions of
cyclic homology and Hochschild homology,
\begin{equation}\label{tw.hoch}
HC_\idot^F(A) = H_\idot(\Lambda,FA_\hash), \qquad HH^F_\idot(A) =
H_\idot(\Delta^o,j^*FA_\hash).
\end{equation}

\begin{exa}\label{cycl.hash}
For the trace functor $\Cycl^l$ of Example~\ref{cycl.exa}, the
functor $\Cycl^lA_\hash:\Lambda \to k\proj$ can be explicitly
described as follows. Let $\Lambda_l$ be the category with the same
objects as $\Lambda_\infty$ and morphisms given by
$$
\Lambda_l([n],[m]) = \Lambda_\infty([n],[m])/\sigma^l,
$$
and let $i_l,\pi_l:\Lambda_l \to \Lambda$ be the natural functors of
\cite[Subsection 1.5]{K} ($\pi_l$ is the natural projection,
$\pi_l([n]) = [n]$, and $i_l$ is the edgewise subdivision functor
sending $[n]$ to $[nl]$). Then we have a natural identification
$$
\Cycl^lA_\hash \cong \pi_{l!}i_l^*A_\hash.
$$
\end{exa}

\subsection{Trace theories.}

Fix a monoidal category $\C$ as in Section~\ref{tr.fun.sec}, with
the additional assumption $0 \in \C$ of
Subsection~\ref{enr.subs}. For any associative unital algebra $A$ in
$\C$, denote by $A\bimod$ the category of $A$-bimodules in $\C$, or
equivalently, the category of left modules over $A \otimes A^o$,
where $A^o$ stands for the opposite algebra. Let $A\biproj \subset
A\bimod$ be the full subcategory of bimodules of the form $V \otimes
A$, where $V$ is a left $A$-module. More generally, for two
associative unital algebras $A$, $B$ in $\C$, denote by
\begin{equation}\label{biproj.eq}
(A,B)\biproj \subset (A \otimes B^o)\amod
\end{equation}
the full subcategory spanned by modules of the form $V \otimes B$,
$V$ a left $A$-module. Explicitly, for any $A$-modules $V$, $V'$, we
have
\begin{equation}\label{A.B.eq}
\Hom_{(A,B)\biproj}(V \otimes B,V' \otimes B) = \Hom_A(V,V' \otimes
B).
\end{equation}
Then for any three associative unital algebras $A$, $B$, $C$ in
$\C$, we have the natural tensor product functor
$$
- \otimes_B -:(A,B)\biproj \times (B,C)\biproj \to (A,C)\biproj
$$
sending $(V \otimes B) \times (W \otimes C)$ to $V \otimes W \otimes
C$. Here maps in $(A,B)\biproj$ act on the product by the
composition
$$
\begin{CD}
\Hom_A(V,V' \otimes B) @>>> \Hom_A(V \otimes W,V' \otimes B \otimes
W) \\@>{m}>> \Hom_A(V \otimes W,V' \otimes W)\\ @>{e}>> \Hom_A(V \otimes
W,V' \otimes W \otimes C),
\end{CD}
$$
where we have used \eqref{A.B.eq}, $m$ is the $B$-module structure
map on $W$, and $e$ is the unity map of the algebra $C$; the action
of the maps in $(B,C)\biproj$ is obvious and left to the reader.
This tensor product functor is equipped with an obvious
associativity isomorphism.

\begin{defn}\label{tr.th.def}
A {\em trace theory} $F$ in $\C$ with values in a category $\E$ is a
collection of functors $F_A:A\biproj \to \E$, one for each
associative unital algebra $A$ in $\C$, and functorial isomorphisms
$$
\tau_{M,N}:F_A(M \otimes_B N) \cong F_B(N \otimes_A M),
$$
one for each pair of associative unital algebras $A$, $B$ and any
bimodules $N \in (A,B)\biproj$, $M \in (B,A)\biproj$, such that for
any three associative unital algebras $A$, $B$, $C$ in $\C$ and any $M \in
(A,B)\biproj$, $N \in (B,C)\biproj$, $L \in (C,A)\biproj$, we have
\eqref{M.N.L}, and for any associative unital algebra $A$ in $\C$
and $M \in A\biproj$, we have $\tau_{A,M} = \tau_{M,A} = \id$.
\end{defn}

In particular, $1 \in \C$ is a unital associative algebra, with
$1\biproj = \C$, and $F_1$ is a functor from $\C$ to $\E$; we say
that a trace theory is {\em normalized} if $\tau_{1,M} = \id$ for
any $M \in 1\biproj = \C$. Then $F_1:\C \to \E$ is a trace functor
on $\C$ with values in $\E$, so that the notion of a trace theory is
a generalization of the notion of a trace functor on $\C$.

However, the generalization is illusory: we have the following.

\begin{prop}\label{tr.fu.th}
Any trace functor $\langle F, \tau_{M,N}\rangle$ from $\C$ to $\E$
uniquely extends to a normalized trace theory in $\C$.
\end{prop}

Before proving this, we need to generalize slightly
Definition~\ref{tw.hash.def}. As in Subsection~\ref{enr.subs}, the
pair of an algebra $A \in \C$ and an $A$-bimodule $M$ corresponds to
a $\C$-enrichment of the category $B(\Nb) = [1]_\Lambda$ of a special
form (effectively, this describes a split square-zero extension of
$A$). More generally, take an object $[n] \in \Lambda$.

\begin{defn}\label{sq.0.enr}
A $\C$-enrichment $\phi$ of the category $[n]_\Lambda$ is {\em
  square-zero} if for any morphism $f$ in $[n]_\Lambda$ of length
$l(f) \geq 2$, we have $\phi(f)=0$.
\end{defn}

Explicitly, such a square-zero enrichment $\phi$ corresponds to a
collection of algebras $A_v = \phi(0_v) \in \C$, one for each $v \in
V([n])$, and a collection of modules $M_v = \phi(1_v) \in (A_v
\otimes A_{v+1}^o)\amod$. Then as in the case $n=1$, $\phi$ induces
a functor
$$
\Lambda(\phi):\Lambda([n]_\Lambda) \to \Lambda(B(\C^o)).
$$
Composing this with the functor $j_n$ of \eqref{j.n}, the natural
embedding $\Lambda/[n] \to \Lambda([n]_\Lambda)$, and the functor
$F_\hash$ of Lemma~\ref{tr.hash}, we obtain a functor from
$(\Delta^o)^n$ to $\E$ which we denote by
\begin{equation}\label{F.multi}
F(M_\idot/A_\idot)_\hash:(\Delta^o)^n \to \E.
\end{equation}

\proof[Proof of Proposition~\ref{tr.fu.th}.] Uniqueness is obvious:
the isomorphisms $F_A(V \otimes A) \cong F(V)$ for any associative
unital algebra $A$ and left $A$-module $V$ are a part of the
definition of a trace theory, and the isomorphisms $\tau_{-,-}$ are
induced by the corresponding isomorphisms for the trace functor
$F$. We have to prove existence.

Let $\Delta_+$ be the category of finite non-empty totally ordered
sets and maps between them sending the initial element to the
initial element. Let $\iota:\Delta_+ \to \Delta$ be the tautological
embedding. Note that the $\iota$ has a left-adjoint functor
$\kappa:\Delta \to \Delta_+$ that adds a new initial element to a
totally ordered set $[n] \in \Delta$. The corresponding functor
$\kappa:\Delta^o \to \Delta^o_+$ between the opposite categories is
then right-adjoint to the tautological functor $\iota:\Delta^o_+ \to
\Delta^o$.

The object $[1] \in \Delta^o_+$ is both terminal and initial. Thus
for any functor $F:\Delta^o_+ \to \E$, the colimit
$$
\lim_{\underset{\Delta_+^o}{\to}}F
$$
exists and is naturally identified with $F([1])$. Then by
adjunction, $\kappa^*$ is adjoint to $\iota^*$, so that we have
$$
\lim_{\underset{\Delta^o}{\to}}\kappa^*F \cong
\lim_{\underset{\Delta^o_+}{\to}}F \cong F([1])
$$
(and in particular, the colimit on the left-hand side exists).

Now note that for any associative unital algebra $A$ in $\C$ and any
$A$-bimodule $M$, the restriction $\iota^*(M/A)_\hash$ of the
functor \eqref{wt.M.A.eq} to the subcategory $\Delta^o_+ \subset
\Delta^o$ only depends on the left $A$-module structure on $M$. In
particular, we can define a natural functor
$$
(V/A)_\hash:\Delta^o_+ \to \C_\hash
$$
for any left $A$-module $V$. Assume given a trace functor $F$ from
$\C$ to $\E$. Then composing the functor $F_\hash$ with
$(V/A)_\hash$, we obtain a natural functor
$$
F(V/A)_\hash = F_\hash \circ (V/A)_\hash:\Delta_+^o \to \E.
$$
Moreover, for $M = V \otimes A$, we have a natural indeitification
$$
F(M/A)_\hash \cong \kappa^*F(V/A)_\hash.
$$
We conclude that
$$
\lim_{\underset{\Delta^o}{\to}}F(M/A)_\hash \cong F(V),
$$
and in particular, the colimit exists.  But the left-hand side
clearly defines a functor from $A\biproj$ to $\E$; this is our
functor $F_A$.

To construct the isomorphism $\tau_{M,N}$ for some $A$, $B$, $M \in
(A,B)\biproj$, $N \in (B,A)\biproj$, let $A_1 = A$, $A_2 = B$, $M_1
= M$, $M_2 = N$, and consider the functor
$$
F(M_\idot/A_\idot)_\hash:\Delta^o \times \Delta^o \to \E
$$
of \eqref{F.multi}. Then one immediately checks that for any object
$[n] \in \Delta^o$, the restriction of this functor to $[n] \times
\Delta^o \subset \Delta^o \times \Delta^o$ extends to the category
$\Delta^o_+ \supset \kappa(\Delta^o)$, and the value of the extended
functor at $[0] \in \Delta^o_+$ is canonically identified with
$F((M \otimes_B N)/A)_\hash([n])$. Therefore the limit
$$
\lim_{\overset{\Delta^o}{\to}}F(M_\idot/A_\idot)_\hash
$$
with respect to the second variable exists and is canonically
identified with $F(M \otimes_B N)_\hash$. We conclude that we have a
natural identification
$$
F_A(M \otimes_B N) = \lim_{\overset{\Delta^o}{\to}}F(M \otimes_B
N)_\hash \cong \lim_{\overset{\Delta^o \times
    \Delta^o}{\longrightarrow}}F(M_\idot/A_\idot)_\hash.
$$
Since the right-hand side is clearly symmetric with respect to
interchanging $A_1$ with $A_2$ and $M_1$ with $M_2$, we obtain the
isomorphism $\tau_{M,N}$. Doing the three-variable version of this
construction, we obtain \eqref{M.N.L}.
\endproof

\begin{remark}
Essentially, for an algebra $A$ and bimodule $M \in A\biproj$, the
object $F_A(M) \in \E$ is defined by the coequalizer diagram
$$
F(M \otimes
A)
\begin{array}{c}\longrightarrow\\[-2mm]\longrightarrow
\end{array}
F(M) \longrightarrow F_A(M),
$$
where the first map on the left-hand side is induced by the product
map $A \otimes M \to M$, and the second map is the composition of
the map $\tau_{A,M}$ and the map induced by the product map $M
\otimes A \to M$. The trickery with the category $\Delta_+^o$ is
needed to show that the coequalizer exists, without any assumptions
on the target category $\E$.
\end{remark}

For the cyclic power trace functor $\Cycl^l$ of
Example~\ref{cycl.exa}, the corresponding trace theory can be
explicitly described as follows. Let $\C$ be the category of flat
$k$-modules, assume given a flat associative unital $k$-algebra $A$
and an $A$-bimodule $M \in A\biproj$, consider the $l$-th tensor
power $A^{\otimes_k l}$, and let $\sigma_A:R^{\otimes_k l} \to
A^{\otimes_k l}$ be the cyclic permutation of order $l$. Then
$A^{\otimes_k l}$ is an associative algebra, $\sigma_A$ is an
algebra map, and we have
\begin{equation}\label{hh.tensor}
HH_0(A,M^{\otimes_A i}) = HH_0(A, M \otimes_A \dots \otimes_A M) =
HH_0(A^{\otimes_k l},M^{\otimes_k l}_\sigma),
\end{equation}
where $M^{\otimes_k l}_\sigma$ is $M^{\otimes_k l}$ as a $k$-module,
with the $R^{\otimes l}$-bimodule structure given by
$$
a_1 \cdot m \cdot a_2 = a_1m\sigma_A(a_2), \qquad a_1,a_2 \in
A^{\otimes_k l},m \in M^{\otimes_k l}.
$$
The cyclic permutation $\sigma_M:M^{\otimes_k l} \to M^{\otimes_k l}$ of
order $l$ together with $\sigma_A$ induce an order-$l$ automorphism
of the right-hand side of \eqref{hh.tensor} which we denote by
$\sigma$.

\begin{lemma}\label{cycl.A}
Let $\Cycl^l_A:A\biproj \to k\proj$ be the trace theory associated
to the functor $\Cycl^l$ by Proposition~\ref{tr.fu.th}. Then we have
$$
\Cycl^l_A(M) = HH_0(A^{\otimes_k l},M^{\otimes_k l}_\sigma)_\sigma.
$$
\end{lemma}

\proof{} A straightforward computation. \endproof

\section{Normalization and denormalization.}\label{N.D.sec}

We now turn our attention to the main subject of the paper, DG
categories. This section contains the necessary preliminaries about
chain complexes.

\subsection{Chain-cochain complexes.}\label{c-c.subs}

Assume given an additive category $\E$ which has countable direct
sums. Denote by $C_\idot(\E)$ the category of (unbounded) complexes
in $\E$. Let $C_\go(\E),C^\go(\E) \subset C_\idot(\E)$ be the full
subcategories spanned by complexes trivial in degrees $<0$
resp. $>0$, where our convention between homological and
cohomological degrees is $M^i=M_{-i}$, $i \in \Z$.

\begin{defn}
The category $C^\go_\go(\E)$ is the category of second-\-quad\-rant
bicomplexes in $\E$, or in other words,
$$
C^\go_\go(\E) = C^\go(C_\go(\E)) \cong C_\go(C^\go(\E)).
$$
\end{defn}

We will call objects in $C^\go_\go(\E)$ {\em chain-cochain
  complexes} in $\E$. Explicitly, chain-cochain complexes are
double-graded objects $M^\hdot_\idot$ in $\E$, with $M^i_j$ being
the term of bidegree $\langle i,j \rangle$, equipped with two
commuting differentials $d_\idot:M_\idot^\hdot \to
M_{\idot-1}^\hdot$, $d^\hdot:M_\idot^\hdot \to M_\idot^{\hdot+1}$
and such that $M^i_j=0$ unless $i,j \geq 0$.

Since $\E$ by assumption has countable direct sums, we have a
natural functor
\begin{equation}\label{tot}
\T:C^\go_\go(\E) \to C_\idot(\E)
\end{equation}
sending a bicomplex $M^\hdot_\idot$ to its sum-total complex
$\T(M)_\idot$ with terms
$$
\T(M)_i =
\begin{cases}
\bigoplus_{j \geq 0}M^j_{j+i}, &\quad i \geq 0,\\
\bigoplus_{j \geq 0}M^{j-i}_j, &\quad i < 0,
\end{cases}
$$
and the differential $d = d^\hdot + (-1)^id_\idot:\T(M)_i \to
\T(M)_{i-1}$.

\begin{defn}\label{ver.equi}
A map $a:M^\hdot_\idot \to N^\hdot_\idot$ between two chain-cochain
complexes in $\E$ is a {\em vertical chain-homotopy equivalence} if
it is a chain-homotopy equivalence with respect to the differential
$d^\hdot$.
\end{defn}

\begin{lemma}\label{ver.con}
If $f:M^\hdot_\idot \to N^\hdot_\idot$ is a vertical chain-homotopy
equivalence, then $\T(f)$ is a chain-homotopy equivalence.
\end{lemma}

\proof{} Let $P^\hdot_\idot$ be the cone of the map $f$ taken in the
cohomological direction, that is,
$$
P^i_j = M^i_j \oplus N^{i-1}_j
$$
with the obvious differentials. Then $P^\hdot_\idot$ is contractible
with respect to the differential $d^\hdot$, and $\T(P^\hdot_\idot)$
is the cone of the map $\T(f)$, so that it suffices to prove that
$\T(P^\hdot_\idot)$ is contractible with respect to its differential
$d=d^\hdot + (-1)^id_\idot$. Indeed, let $h:P^{\hdot+1}_\idot \to
P^\hdot_\idot$ be a map such that $hd^\hdot + d^\hdot h = \id$. Then
$\nu=\id-(\T(h)d+d\T(h)):\T(P^\hdot_\idot) \to \T(P^\hdot_\idot)$
decreases the homological degree; therefore the series
$$
p = \id + \nu + \dots + \nu^n + \dots
$$
is well-defined on $\T(P^\hdot_\idot)$. We then have $d(ph) + (ph)d
= p(dh+hd)=p(\id-\nu) = \id$, so that $ph$ is a contracting homotopy
for $\T(P^\hdot_\idot)$.
\endproof

The functor $\T$ has a right-adjoint
$$
\I:C_\idot(\E) \to C^\go_\go(\E);
$$
explicitly, if a complex $M_\idot \in C_\idot(\E)$ consists of a
single object $M$ placed in degree $i$, then $\I(M_\idot)$ is the
bicomplex
\begin{equation}\label{I.exp}
\begin{CD}
\dots\\
@AA{\id}A\\
M @>{\id}>> M\\
@. @AA{\id}A\\
@. M @>{\id}>> M\\
@. @. @AA{\id}A\\
@. @. M @>{\id}>> M,
\end{CD}
\end{equation}
with $M$ in the bottom right corner placed in bidegree $\langle i,0
\rangle$ if $i \geq 0$ and $\langle 0,-i \rangle$ otherwise. 

\begin{lemma}\label{tot.le}
For any complex $M_\idot \in C_\idot(\E)$ in the category $\E$, the
adjunction map $\T(\I(M_\idot)) \to M_\idot$ is a split surjection
and a chain-homotopy equivalence.
\end{lemma}

\proof{} To obtain a splitting, consider another functor
$\iota:\C_\idot(\E) \to C^\go_\go(\E)$ which only keeps the bottom
right corner of \eqref{I.exp}: for any $M_\idot \in C_\idot(\E)$,
$\iota(M_\idot)$ is a bicomplex $M^\hdot_\idot$ with $M_i^0 = M_i$,
$M_0^j = M_{-j}$, $M^i_j=0$ if $i,j \geq 1$, and the obvious
differentials. We have a natural isomorphism $\T \circ \iota \cong
\id$, and by adjunction, this induces a map $f:\iota \to
\I$. Applying $\T$, we obtain a map
$$
\T(f):\id \cong \T \circ \iota \to \T \circ \I
$$
left-inverse to the adjunction map $\T \circ \I \to \id$. 

To prove that the adjunction map is a chain-homotopy equivalence, it
then suffices to prove that $\T(f):M_\idot \to \T(\I(M_\idot))$ is a
chain-homotopy equivalence, and by Lemma~\ref{ver.con}, it suffices
to prove that $f:\iota(M_\idot) \to \I(M_\idot)$ is a vertical
chain-homotopy equivalence.

Assume first that $\E = \Z\proj$ is the category of projective
$\Z$-modules, and let $M_\idot$ be $\Z$ placed in degree $0$. Then
by \eqref{I.exp},
\begin{equation}\label{I-de}
I = \I(\Z) \cong \Z[u]\langle\eps\rangle
\end{equation}
is the free commutative bi-DG algebra generated by $u$ in bidegree
$\langle 1,1 \rangle$ and $\eps$ in bidegree $\langle 0,1 \rangle$,
with $d^\hdot(\eps)=u$ and $d_\idot(\eps)=1$. On the other hand,
$\iota(\Z)$ is $\Z$ placed in bidegree $\langle 0,0 \rangle$, and
the map $f:\iota(\Z) \to \I(\Z)$ is obviously a vertical
chain-homotopy equivalence. Now for general $\E$, $M_\idot$, it is
straightforward to see that we have an isomorphism
\begin{equation}\label{I-I}
\I(M_\idot) \cong \iota(M_\idot) \otimes I,
\end{equation}
and the statement for $\Z$ implies the statement for $M_\idot$.
\endproof

\begin{defn}\label{hor-ver}
Two maps $f_1,f_2:N^\hdot_\idot \to M^\hdot_\idot$ between
chain-cochain complexes $N^\hdot_\idot,M^\hdot_\idot \in
C^\go_\go(\E)$ are {\em horizontally chain-homotopic} if the exists
a map $h:N^\hdot_\idot \to M^\hdot_{\idot+1}$ such that
$$
h d^\hdot + d^\hdot h=0, \qquad h d_\idot + d_\idot h = f_1-f_2.
$$
A map $f:N^\hdot_\idot \to M^\hdot_\idot$ is a {\em horizontal
chain-homotopy equivalence} if the exists a map $f':M^\hdot_\idot
\to N^\hdot_\idot$ such that $f \circ f'$ and $f' \circ f$ are
horizontally chain-homotopic to the identity maps.
\end{defn}

\begin{lemma}\label{hori.homo}
Two maps $f_1,f_2:N_\idot \to M_\idot$ between two complexes
$N_\idot,M_\idot \in C_\idot(\E)$ are chain-homotopic if and only if
$\I(f_1)$ and $\I(f_2)$ are horizontally chain-homotopic in the
sense of Definition~\ref{hor-ver}.
\end{lemma}

\proof{} The ``if'' part is obvious: a horizontal chain homotopy
between $\I(f_1)$ and $\I(f_2)$ induces a chain homotopy between
$\T(\I(f_1))$ and $\T(\I(f_2))$, and this is enough by
Lemma~\ref{tot.le}.

For the ``only if'' part, note that since the functor $\I$ is
additive, it suffices to prove that if $f=f_1-f_2$ is
chain-homotopic to $0$, then $\I(f)$ if horizontally chain-homotopic
to $0$. For any complex $E_\idot \in C_\idot(\E)$, let $C(E_\idot)$
be the cone of the identity endomorphism of $E_\idot$ -- that is,
the natural extension
$$
\begin{CD}
0 @>>> E_\idot @>>> C(E_\idot) @>>> E_\idot[1] @>>> 0.
\end{CD}
$$
Then $C(E_\idot)$ is canonically contractible for any $E_\idot$, in
particular for $E_\idot = N_\idot$, and $f:N_\idot \to M_\idot$ is
homotopic to $0$ if and only if it factors through the embedding
$N_\idot \to C(N_\idot)$. On the other hand, $C$ sends $C_\go(\E)$
to itself, thus induces an endofunctor of $C_\go(\E)$ and of
$$
C^\go_\go(\E) = C_\go(C^\go(\E)).
$$
Denote this endofunctor by $C_{<}$. Then $\I(f)$ is horizontally
chain-homotopic to $0$ if and only if it factors through the natural
embedding $c:\I(N_\idot) \to C_{<}(\I(N_\idot))$. By adjunction, we
have to show that there exists a map $\wt{f}:\T(C_{<}(\I(N_\idot)))
\to M_\idot$ such that the diagram
$$
\begin{CD}
\T(\I(N_\idot)) @>{\ad}>> N_\idot\\
@V{\T(c)}VV @VV{f}V\\
\T(C_{<}(\I(N_\idot))) @>{\wt{f}}>> M_\idot
\end{CD}
$$
is commutative. However, we have an obvious isomorphism $\T \circ
C_{<} \cong C \circ \T$. Thus if $f$ is chain-homotopic to $0$, we
can define $\wt{f}$ as the composition
$$
\begin{CD}
\T(C_{<}(\I(N_\idot))) \cong C(\T(\I(N_\idot))) @>{C(\ad)}>>
C(N_\idot) \to M_\idot,
\end{CD}
$$
where the right-hand side map extends $f$.
\endproof

Assume now for a moment that $\E$ is a symmetric tensor
category. Then all the categories $C_\idot(\E)$, $C_\go(\E)$,
$C^\go(\E)$, $C_\go^\go(\E)$ are also symmetric tensor categories,
and the total complex functor $\T:C^\go(k) \to C_\idot(k)$ of
\eqref{tot} is a tensor functor. By adjunction, the functor
$\I:C_\idot(k) \to C^\go_\go(k)$ is pseudotensor -- namely, for any
$M_\idot,N_\idot \in C_\idot(k)$, the map
$$
\begin{CD}
\T(\I(M_\idot) \otimes \I(N_\idot)) \cong \T(\I(M_\idot)) \otimes
\T(\I(N_\idot)) @>{\ad \otimes \ad}>> M_\idot \otimes N_\idot
\end{CD}
$$
induces a functorial map
\begin{equation}\label{I.ten}
\I(M_\idot) \otimes \I(N_\idot) \to \I(M_\idot \otimes N_\idot)
\end{equation}
satisfying the usual compatibilities. This
map becomes a chain-homotopy equivalence after applying the total
complex functor $\T$.

\subsection{Dold-Kan equivalences.}\label{D.K.subs}

Now assume that the additive category $\E$ is Karoubi-closed, and
let $\Delta$, $\Delta^o$ be category of non-empty finite totally
ordered sets and its opposite, as in
Subsection~\ref{tr.def.subs}. Recall that for any $M \in \Delta^o\E$
the {\em normalized chain complex} $\Nn(M)$ is a certain canonical
subcomplex in the standard complex $M_\idot$, and sending $M$ to
$\Nn(M)$ induces an equivalence of categories
$$
\Nn_\idot:\Delta^o\E \to C_\go(E),
$$
the Dold-Kan equivalence (see e.g. \cite[Section III.2]{gj} for a
proof). Replacing $\E$ with the opposite category $\E^o$ and using
the obvious identifications
\begin{equation}\label{opp.id}
(\Delta^o\E^o)^o = \Delta\E, \qquad C_\go(E^o)^o = C^\go(\E),
\end{equation}
we obtain an equivalence
$$
\Nn^\hdot:\Delta\E \to C^\go(E)
$$
for cosimplicial objects. Taking both constructions together, we
obtain an equivalence
$$
\Nn:(\Delta^o \times \Delta)\E \to C^\go_\go(E),
$$
where $(\Delta^o \times \Delta)\E$ is the category of functors from
$\Delta^o \times \Delta$ to $\E$. Denote by
$$
\Dd_\idot:C_\go(\E) \to \Delta^o\E, \quad \Dd^\hdot:C^\go(\E) \to
\Delta\E, \quad \Dd:\C^\go_\go(E) \to (\Delta^o \times \Delta)\E
$$
the denormalization functors giving inverse equivalences. Combining
this with the functors $\T$ and $\I$, we obtain a pair of adjoint
functors
\begin{equation}\label{tot.nd}
\Ll=\T \circ \Nn:(\Delta^o \times \Delta)\E \to C_\idot(E), \quad
\Rr=\Dd \circ \I:\C_\idot(E) \to (\Delta^o \times \Delta)\E,
\end{equation}
and by Lemma~\ref{tot.le}, the adjunction map
$$
\Ll(\Rr(M_\idot)) \to M_\idot
$$
is a chain-homotopy equivalence for any $M_\idot \in C_\idot(\E)$.

The normalized chain complex $\Nn_\idot(M)$ of a simplicial object
$M \in \Delta^o\E$ is a chain-homotopy retract of the standard
complex $M_\idot$ -- that is, $\Nn_\idot(M)$ is a direct summand in
$M_\idot$, and the embedding $\Nn_\idot(M) \to M_\idot$ is a
chain-homotopy equivalence. Analogously, for any $M \in \Delta\E$,
$\Nn^\hdot(M)$ is a chain-homotopy retract of the standard complex
$M^\hdot \in C^\go(\E)$.

If $\Delta_+ \subset \Delta$ is the subcategory of
Proposition~\ref{tr.fu.th}, then the Dold-Kan equivalence induces an
equivalence between $\Delta^o_+\E$ and the category of objects in
$\E$ graded by non-negative integers, and restricting to
$\Delta^o_+$ corresponds to forgetting the grading (it can be proved,
for example, by exactly the same argument as in \cite[Section
  III.2]{gj}, namely, deduced from \cite[Proposition
  III.2.2]{gj}). Analogously, $\Delta\E$ is also equivalent to the
category of graded objects, and restricting to $\Delta_+ \subset
\Delta$ corresponds to forgetting the grading. In fact, $\Delta_+$
is equivalent to $\Delta_+^o$.

If we extend the normalization functor $\Nn_\idot$ to a functor
\begin{equation}\label{double.n}
\Nn^{[2]}_\idot:(\Delta^o \times \Delta^o)\E \to C_\go(\E)
\end{equation}
by applying normalization in each of the two simplicial directions
and then taking the total complex of the resulting bicomplex, then
for any $B \in (\Delta^o \times \Delta^o)\E$, we have natural maps
\begin{equation}\label{big.shuffle}
\sh_\idot:\Nn_\idot(\delta^*B) \to \Nn^{[2]}_\idot(B), \quad
\sh'_\idot:\Nn^{[2]}_\idot(B) \to \Nn_\idot(\delta^*B),
\end{equation}
where $\delta:\Delta^o \to \Delta^o \times \Delta^o$ is the diagonal
embedding. The maps $\sh_\idot$, $\sh'_\idot$ are known as the {\em
  Eilenberg-Zilber} and the {\em Alexander-Whitney shuffle}. They
are associative in the obvious sense, and they are mutually inverse
chain-homotopy equivalences. 

\bigskip

Assume now given two additive Karoubi-closed categories $\C$, $\E$,
and a functor $F:\C \to \E$. Applying $F$ pointwise, we extend it to
a functor $F^{\Delta}:\Delta\C \to \Delta\E$, and we let
$$
F^\hdot = \Nn^\hdot \circ F^\Delta \circ \Dd^\hdot:C^\go(\C) \to
C^\go(\E)
$$
be the corresponding functor between the categories of
complexes. Then the famous classic theorem of Dold \cite{dold}
claims that $F^\hdot$ sends chain-homotopic maps to chain-homotopic
maps.  Passing to the opposite categories $\C^o$, $\E^o$, we get the
same statement for $F^{\Delta^o}:\Delta^o\C \to \Delta^o\E$ and
$F_\idot = \Nn_\idot \circ F^{\Delta^o} \circ \Dd_\idot$.  In
particular, $F^\hdot$ and $F_\idot$ preserve chain-homotopy
equivalences. All of this requires no assumptions on $F$ whatsoever.

For chain-cochain complexes, this has the following
corollary. Extend $F$ to a functor $F^{\Delta^o \times
  \Delta}:(\Delta^o \times \Delta)\C \to (\Delta^o \times
\Delta)\E$, again by applying it pointwise, and let
$$
F^\hdot_\idot = \Nn \circ F^{\Delta^o \times \Delta} \circ
\Dd:C^\go_\go(\C) \to C^\go_\go(\E).
$$

\begin{lemma}\label{ho-ve}
Assume that a map $f:M^\hdot_\idot \to N^\hdot_\idot$ between
chain-cochain complexes is either a vertical chain-homotopy
equivalence in the sense of Definition~\ref{ver.equi}, or a
horizontal chain-homotopy equivalence in the sense of
Definition~\ref{hor-ver}. Then the map $F^\hdot_\idot(f)$ is also a
vertical resp. a horizontal chain-homotopy equivalence.
\end{lemma}

\proof{} In the horizontal case, this immediately follows from
applying Dold Theorem in the simplicial direction. In the vertical
case, note that a map is a vertical chain-homotopy equivalence if
and only if it is becomes a chain-homotopy equivalence after we
forget the horizontal differential. Under the Dold-Kan equivalence,
forgetting the horizontal differential correponds to restricing to
the subcategory $\Delta_+^o \times \Delta$; after that, we can apply
Dold Theorem in the cosimplicial direction and get the claim.
\endproof

\subsection{Balanced functors.}

To pass from chain-cochain complexes to unbounded complexes, assume
that $\C$ and $\E$ have countable sums, and moreover, assume that
$\E$ is equipped with a structure of an exact category.

\begin{defn}\label{tot.def}
A map $f:M^\hdot_\idot \to N^\hdot_\idot$ between two chain-cochain
complexes is a {\em total equivalence} if $\T(f)$ is a
chain-homotopy equivalence.
\end{defn}

\begin{defn}\label{adm.fun.def}
A functor $F:\C \to \E$ is {\em balanced} if for any total
equivalence $f:M^\hdot_\idot \to N^\hdot_\idot$,
$M^\hdot_\idot,N^\hdot_\idot \in \C^\go_\go\C$, the map
$\T(F^\hdot_\idot(f))$ is a quasiisomorphism.
\end{defn}

Unfortunately, it seems that this definition is not vacuous -- there
is no analog of a Dold Theorem for unbounded chain complexes and
arbitrary functors $F$. Therefore we need some results on balanced
functors.

\begin{lemma}\label{ba.iso}
A functor $F:\C \to \E$ is balanced if and only if for any
chain-cochain complex $M^\hdot_\idot \in C_\go^\go(\C)$ with
$M^\hdot = \T(M^\hdot_\idot)$ and the adjunction map
$a:M^\hdot_\idot \to \I(M^\hdot)$, the induced map
$$
\T(F^\hdot_\idot(a)):\T(F^\hdot_\idot(M^\hdot_\idot)) \to
\T(F^\hdot_\idot(\I(M^\hdot)))
$$
is a quasiisomorphism.
\end{lemma}

\proof{} One direction is obvious -- by Lemma~\ref{tot.le}, $a$ is a
total equivalence in the sense of Definition~\ref{tot.def}. In the
other direction, by Definition~\ref{adm.fun.def}, it suffices to
show that $\Ll(F(\Rr(\T(f))))$ is a quasiisomorphism for an total
equivalence $f:M^\hdot_\idot \to N^\hdot_\idot$; this follows by
applying Lemma~\ref{hori.homo}, Lemma~\ref{ho-ve}, and then
Lemma~\ref{ver.con} and Lemma~\ref{hori.homo} again.
\endproof

\begin{lemma}
A pointwise extension and a pointwise filtered direct limit of
balanced functors is balanced.
\end{lemma}

\proof{} Clear. \endproof

Thus balanced functors can be built out of building blocks. One
class of such blocks is provided by the following technical
gadget. For any additive category $\C$, let $W_\idot\C$ be the
category of objects in $\C$ equipped with a finite-length increasing
term-wise split filtration numbered by non-negative integers, and
filtered maps between these objects. In other words, objects of
$W_\idot\C$ are collections $C_\idot=\langle C_0,C_1,\dots,C_n,\dots
\rangle$ of objects in $\C$, and maps are given by
$$
\Hom(C_\idot,C'_\idot) = \prod_{i \geq j \geq 0}\Hom(C_i,C'_j).
$$
We have an obvious forgetful functor $\tau:W_\idot\C \to \C$ sending
$C_\idot$ to $\bigoplus_i C_i$, and acting on morphisms by the
natural embedding
$$
\prod_{i \geq j \geq 0}\Hom(C_i,C'_j) \subset
\Hom\left(\bigoplus_i C_i,\bigoplus_j C_j\right).
$$
We also have the associated graded object functor $\gr:W_\idot\C \to
\C$, again sending $C_\idot$ to $\bigoplus_iC_i$ and acting on
morphisms by the natural map
$$
\prod_{i \geq j \geq 0}\Hom(C_i,C'_j) \to
\prod_i\Hom(C_i,C_i) \subset \Hom\left(\bigoplus_i
C_i,\bigoplus_j C_j\right).
$$

\begin{defn}\label{bounded.defn}
A functor $F:\C \to \E$ between additive categories $\C$ and $\E$ is
{\em filterable} if there exists a functor $W_\idot F:W_\idot\C \to
W_\idot\E$ such that $F \circ \gr \cong \gr \circ W_\idot F$ and $F
\circ \tau \cong \tau \circ W_\idot F$.
\end{defn}

Informally speaking, a split filtration on an object $C \in \C$
should functorially induce a split filtration on $F(C) \in \E$ such
that $\gr F(C) \cong F(\gr C)$.

\begin{prop}\label{fil=>bal}
Assume that the functor $F:\C \to \E$ is filterable in the sense of
Definition~\ref{bounded.defn}. Then it is balanced in the sense of
Definition~\ref{adm.fun.def}.
\end{prop}

\proof{} Take some $M^\hdot_\idot \in C^\go_\go(C)$ with $M^\hdot =
\T(M^\hdot_\idot)$. Let $I = \Z[u]\langle\eps\rangle$ be as in the
proof of Lemma~\ref{tot.le}, so that $\I(M^\hdot_\idot) \cong I
\otimes \iota(M^\hdot)$. Extend the adjunction map to an $I$-module
map
$$
\psi:I \otimes M^\hdot_\idot \to I \otimes \iota(M^\hdot) \cong
\I(M^\hdot).
$$
Note that this map is injective, and it becomes an isomorphism after
we invert the generator $u \in I$. Moreover, for any $i \geq 0$, let
$$
W_i\I(M^\hdot) = u^{-i}\psi(I \otimes M^\hdot_\idot) \subset
\I(M^\hdot).
$$
Then this defines an increasing termwise-split filtration on
$\I(M^\hdot)$, and for any $i \geq 1$, the associated graded
quotient $\gr_i^W$ of this filtration is a free module over $I/uI
\cong \Z\langle \eps \rangle$. Since $I/uI$ is obviously
horizontally contractible, $\gr_i^W$ is also horizontally
contractible, so that the map
$$
\gr(\psi):I \otimes M^\hdot_\idot \to \gr^\hdot_W\I(M^\hdot)
$$
is a horizontal chain-homotopy equivalence.

Now consider the filtered extension $W_\idot F$ of the functor $F$
provided by Definition~\ref{bounded.defn}. Then by
Lemma~\ref{ho-ve},
$$
\gr((W_\idot F)^\hdot_\idot(\psi)) \cong F^\hdot_\idot(\gr(\psi))
$$
is also a horizontal chain-homotopy equivalence, and this implies
that the chain-cochain complex $F^\hdot_\idot\I(M^\hdot)$ has an
increaing filtration $W_\idot$ such that $W_0$ is identified with
$F^\hdot_\idot(I \otimes M^\hdot_\idot)$ by the map $\psi$, and
$\gr_i^W$ is horizontally contractible for any $i \geq 1$. Applying
the total complex functor functor $\T$, we conclude that
$\T(F^\hdot_\idot(\I(M^\hdot)))$ carries a filtration
$W_\idot$ such that $\gr_i^W$ is contractible for $i \geq 1$, and
then the map
$$
\T(F^\hdot_\idot(\psi)):\T(F^\hdot_\idot(I \otimes M^\hdot_\idot))
\to \T(F^\hdot_\idot(\I(M^\hdot)))
$$
is a chain-homotopy equivalence by the same argument as in the proof
of Lemma~\ref{tot.le}. It remains to recall that the embedding
$M^\hdot_\idot \to I \otimes M^\hdot_\idot$ is a vertical
chain-homotopy equivalence, so that by Lemma~\ref{ho-ve}, the
map
$$
\T(F^\hdot_\idot(M^\hdot_\idot)) \to \T(F^\hdot_\idot(I \otimes
M^\hdot_\idot))
$$
is a chain-homotopy equivalence, too.
\endproof

\subsection{Shuffle products.}\label{shuffle.subs}

Fix a commutative ring $k$, and consider the category $k\proj$ of
projective $k$-modules. It is a Karoubi-closed additive category
with a symmetric tensor product. To simplify notation, let
$$
\begin{aligned}
C_\idot(k)&=C_\idot(k\proj), \quad C_\go(k)=C_\go(k\proj),\\
C^\go(k)&=C^\go(k\proj), \quad
C^\go_\go(k)=C^\go_\go(k\proj).
\end{aligned}
$$
The categories $\Delta^o k\proj$, $\Delta k\proj$, $(\Delta^o \times
\Delta)k\proj$ are symmetric tensor categories with respect to
pointwise tensor product. The normalization and denormalization
functors of \eqref{D.K.subs} are not tensor. However, they are
pseudotensor, with the pseudotensor structure induced by shuffle
maps \eqref{big.shuffle}. Namely, taking $\E = k\proj$ and $B = M
\boxtimes N$ for some $M,N \in \Delta^ok\proj$, we obtain a
chain-homotopy equivalence
$$
\Nn_\idot(M \otimes N) \cong \Nn_\idot(\delta^*B) \to
\Nn^{[2]}_\idot(B) \cong \Nn_\idot(M) \otimes \Nn_\idot(N).
$$
This map is functorial in $M$ and $N$, and defines a pseudotensor
structure on the normalization functor $\Nn_\idot$. Replacing $\E$
with the opposite category and dualizing, we obtain a pseudotensor
structure on the normalization functor $\Nn^\hdot$. By adjunction,
we then obtain functorial associative shuffle maps
\begin{equation}\label{D.ten.dot}
\begin{aligned}
\sh_\idot:\Dd_\idot(-) \otimes \Dd_\idot(-) &\to \Dd_\idot( -
\otimes -),\\
\sh^\hdot:\Dd^\hdot(-) \otimes \Dd^\hdot(-) &\to \Dd^\hdot( -
\otimes -),
\end{aligned}
\end{equation}
and their tensor product, a functorial associative shuffle map
\begin{equation}\label{D.ten}
\sh:\Dd(M^\hdot_\idot) \otimes \Dd(N^\hdot_\idot) \to
\Dd(M^\hdot_\idot \otimes N^\hdot_\idot).
\end{equation}
The maps $\sh_\idot$, $\sh^\hdot$ of \eqref{D.ten.dot} become
chain-homotopy equivalences after applying normalization, so that
the map $\sh$ of \eqref{D.ten} is a total equivalence in the sense
of Definition~\ref{tot.def}.

\bigskip

Although we will not need this in this paper, let us describe the
map $\sh$ explicitly, for the sake of completeness and for the
convenience of the reader.

\bigskip

For any object $[n]^o \times [m] \in \Delta^o \times \Delta$, let
$k_{n,m} \in (\Delta^o \times \Delta)k\proj$ be the corresponding
representable functor,
$$
k_{n,m}([n']^o \times [m']) = k\left[\Hom_{\Delta^o \times \Delta}([n]^o
  \times [m],[n']^o \times [m])\right],
$$
where for any $[l] \in \Delta$, $[l]^o$ stands for the same object
considered as an object of the opposite category $\Delta^o$. Let
\begin{equation}\label{Q.eq}
Q^\hdot_\idot(n,m) = \Nn(k_{n,m})
\end{equation}
be its normalization. We can treat $Q^\hdot_\idot(-,-)$ as a functor
in $(\Delta^o \times \Delta)C^\go_\go(k)$; this is the kernel of the
normalization functor $\Nn$, in the sense that we have adjunction
isomorphisms
$$
\Dd(N^\hdot_\idot)([n]^o \times [m]) =
\Hom_{C^\go_\go(k)}(Q^\hdot_\idot(n,m),N^\hdot_\idot)
$$
for any $N^\hdot_\idot \in C^\go_\go(k)$. For any
$N^\hdot_\idot,M^\hdot_\idot \in C^\go_\go(k)$, we then have a natural
map
$$
\begin{aligned}
(\Dd(N^\hdot_\idot) \otimes \Dd(M^\hdot_\idot))&([n]^o \times [m]) =
\Dd(N^\hdot_\idot)([n]^o \times [m]) \otimes
\Dd(M^\hdot_\idot)([n]^o \times [m]) \\
&\cong \Hom_{C^\go_\go}(Q^\hdot_\idot(n,m),N^\hdot_\idot) \otimes
\Hom_{C^\go_\go}(Q^\hdot_\idot(n,m),M^\hdot_\idot) \\
&\to \Hom_{C^\go_\go}(Q^\hdot_\idot(n,m) \otimes
Q^\hdot_\idot(n,m),N^\hdot_\idot \otimes M^\hdot_\idot),
\end{aligned}
$$
and to obtain the shuffle map $\sh$ of \eqref{D.ten}, it suffices to
construct maps
$$
Q^\hdot_\idot(n,m) \to Q^\hdot_\idot(n,m) \otimes Q^\hdot_\idot(n,m)
$$
turning $Q^\hdot_\idot(n,m)$ into a coassociative coalgebra object
in $C^\go_\go(k)$. These maps should be functorial in $[n]^o \times
[m]$. Equivalently, one can construct a functorial associative
algebra structure on the dual chain-cochain complexes
$P^\hdot_\idot(n,m)$ given by
\begin{equation}\label{P.eq}
P^i_j(n,m) = \Hom_k(Q^j_i(n,m),k).
\end{equation}
Moreover, let
$$
P^i(n) = \Hom_k(Q_i(n),k), \qquad P_j(m) = \Hom_k(Q^j(m),k),
$$
where $Q_\idot(n) = \Nn_\idot(k_n)$, $Q^\hdot = \Nn^\hdot(k_m)$ are
the normalizations of the representable functors $k_n \in \Delta^o
k\proj$, $k_m \in \Delta k\proj$. Then
$$
P^\hdot_\idot(n,m) \cong P^\hdot(n) \otimes P_\idot(m),
$$
and it suffices to construct DG algebra structures separately on
$P^\hdot(n)$ and $P_\idot(m)$, $[n]^o \in \Delta^o$, $[m] \in
\Delta$.

For $P^\hdot(n)$, we note that as a complex, $P^\hdot(n)$ is by
definition the normalized chain complex of an elementary
$(n-1)$-simplex. In particular, $P^\hdot(1)$ is $k$ placed in
degree $1$. This has an obvious DG algebra structure. For $n \geq
2$, assume by induction that a DG algebra structure on
$P^\hdot(n-1)$ is already constructed. Then as a complex,
$P^\hdot(n)$ fits into a short exact sequence
$$
\begin{CD}
0 @>>> P^\hdot(n-1)[-1] @>>> P^\hdot(n) @>>> P^\hdot(n-1) \oplus k
@>>> 0,
\end{CD}
$$
with the connecting differential $k \oplus P^\hdot(n-1) \to
P^\hdot(n-1)$ being the sum of the unity embedding $k \to
P^\hdot(n-1)$ and the identity map. Then as an algebra,
$P^\hdot(n)$ is a trivial square-zero extension of $k \oplus
P^\hdot(n-1)$ by the bimodule $P^\hdot(n-1)[-1]$, where
$P^\hdot(n-1)$ acts on $P^\hdot(n-1)[-1]$ trivially on the left and
tautologically on the right, while $k$ acts trivially on the right
and by scalar multiplication on the left.

This inductive description is somewhat unwieldy; an alternative
description is provided by the following result.

\begin{lemma}
The category of DG modules over the DG algebra $P^\hdot(n)$ is
equivalent to the category of complex $M_\idot$ of $k$-modules
equipped with a grading
$$
M_\idot = \bigoplus_{1 \leq i \leq n}M^i_\idot
$$
such that the differential $d:M_\idot \to M_{\idot-1}$ sends
$M^i_\idot$ into the sum of $M^j_\idot$ with $j \geq i$.
\end{lemma}

\proof{} Immediately follows by induction.
\endproof

For $P_\idot(m)$, we note that
\begin{equation}\label{P.i}
P_\idot(m) = P_\idot(1)^{\otimes m},
\end{equation}
and $P_\idot(k) = C(k)$ is the cone complex: we have $P_i(k) = k$
for $i=0,1$ and $0$ otherwise, and the differential $d:k \to k$ is
the idenitity map. This has an obvious DG algebra structure; the one
on $P_\idot(m)$ is then induced by \eqref{P.i}.

\begin{remark}
We note that the DG algebras $P_\idot(m)$ are commutative. This
implies that the map $\sh^\hdot$ of \eqref{D.ten.dot} -- that is,
the cosimplicial part of the shuffle map \eqref{D.ten} -- is not
only associative but also commutative. For the map $\sh_\idot$, this
is not true.
\end{remark}

\section{DG algebras.}\label{dg.alg.sec}

\subsection{Twisted homology for DG algebras.}

As in Subsection~\ref{shuffle.subs}, fix a commutative ring $k$, and
consider the category $C_\idot(k)$ of complexes of projective
$k$-modules. This is a unital symmetric tensor category. By a {\em
  DG algebra} $A_\idot$ we will understand a unital associative
algebra object in $C^\hdot(k)$ -- in other words, we only consider
DG algebras over $k$ that are projective as complexes of
$k$-modules. Analogously, we will only consider DG modules and DG
bimodules that are projective as complexes of $k$-modules.

By a {\em $\dd$-algebra} $\A$ we will understand a unital
associative algebra object in the category $(\Delta^o \otimes
\Delta)k\proj$ -- or equivalently, a functor from $\Delta^o \times
\Delta$ to the category of unital associative $k$-algebras that
are projective as $k$-modules. Modules and bimodules over
$\dd$-algebras are also assumed to be pointwise projective as
$k$-modules.

As in Remark~\ref{alg.exa}, in terms of Definition~\ref{enr.def},
DG algebras correspond to $C_\idot(k)$-enrichments of the point
category $\ppt$, and $\dd$-algebras correspond to its $(\Delta^o
\times \Delta)k\proj$-enrichments. Morphisms between DG algebras
resp. $\dd$-algebras correspond to enriched functors from $\ppt$ to
itself (such a functor is of course trivial, but its enrichment
encodes the morphism).

To handle bimodules, we proceed as in Subsection~\ref{enr.subs} --
the pair $\langle A_\idot,M_\idot \rangle$ corresponds to a
square-zero $C_\idot(k)$-enrichment of the category $[1]_\Lambda$,
and analogously for $\dd$-bimodules. We note that both categories
$C_\idot(k)$ and $(\Delta^o \times \Delta)k\proj$ do have a
terminal object $0$, so that square-zero enrichments are
well-defined. $C_\idot(k)$-Enrichments of the identity functor
$[1]_\Lambda \to [1]_\Lambda$ correspond to morphisms of pairs
$\langle A_\idot,M_\idot \rangle \to \langle A'_\idot,M'_\idot
\rangle$ -- that is, to pairs of a DG morphism $f_A:A_\idot \to
A'_\idot$ and a $k$-linear map of complexes $f_M:M_\idot \to
M'_\idot$ such that $f(am)=f(a)f(m)$ and $f(ma)=f(m)f(a)$ for any $a
\in A_\idot$, $m \in M_\idot$. We will say that $f$ is a {\em
  quasiisomorphism} if so are both its components $f_A$ and
$f_M$. Analogously, $(\Delta^o \times \Delta)k\proj$-enrichments of
$\id:[1]_\Lambda \to [1]_\Lambda$ correspond to morphisms of pairs
$f:\langle \A,\M \rangle \to \langle \A',\M' \rangle$; we will say
that $f$ is a {\em total equivalence} if both its components are
total equivalences in the sense of Definition~\ref{tot.def}.

The pseudotensor functors $\Rr$, $\Ll$ of \eqref{tot.nd} allow to
pass between DG algebras and $\dd$-algebras -- for any DG algebra
$A_\idot$ corresponding to an enrichment $\alpha$, we have a natural
$\dd$-algebra $\A = \Rr(A_\idot)$ corresponding to the enrichment
$\Rr \circ \alpha$, and conversely, for any $\dd$-algebra $\A$, we
have a natural DG algebra $A_\idot = \Ll(\A)$. The same works for
pairs of an algebra and a bimodule.

\medskip

Assume now given an exact category $\E$ with countable filtered
colimits, with derived category $\D(\E)$, and assume further that we
are given a trace functor $F:k\proj \to \E$ in the sense of
Definition~\ref{tr.fu.def}. As in Subsection~\ref{D.K.subs}, extend
$F$ to a functor
$$
F^{\Delta^o \times \Delta}:(\Delta^o \times \Delta)k\proj \to
(\Delta^o \times \Delta)\E
$$
by applying it pointwise. Then $F^{\Delta^o \times \Delta}$ is also
a trace functor. Thus for any small $\dd$-algebra $\A$, we have the
natural functor
$$
F^{\Delta^o \times \Delta}A_\hash:\Lambda \to (\Delta^o \times
\Delta)\E
$$
of Definiton~\ref{tw.hash.def}. Composing it with the functor
$\Ll:(\Delta^o \times \Delta)\E \to C_\idot(\E)$ of \eqref{tot.nd},
we obtain a functor
\begin{equation}\label{FA.dd}
F_\idot\A_\hash = \Ll \circ F^{\Delta^o \times \Delta}A_\hash:\Lambda
\to C_\idot(\E).
\end{equation}
Analogously, assume given a bimodule $\M$ over $\A$, and consider
the natural functor
\begin{equation}\label{Fdd.mu}
F^{\Delta^o \times \Delta}\mu_\hash:\Delta^o \to (\Delta^o \times
\Delta)\E
\end{equation}
of Definition~\ref{tw.hash.def}. Then we can compose it with
$\Ll$ and obtain a functor
\begin{equation}\label{FAM.dd}
F_\idot(\M/\A)_\hash = \Ll \circ F^{\Delta^o \times
  \Delta}(\M/\A)_\hash:\Delta^o \to C_\idot(\E).
\end{equation}

\begin{defn}
\begin{enumerate}
\item The {\em twisted cyclic homology object} $CC^F_\idot(\A) \in
  \D(\E)$ of a small $\dd$-algebra $\A$ is given by
$$
CC^F_\idot(\A) = H_\idot(\Lambda,F_\idot\A_\hash),
$$
where $F_\idot\A_\hash$ is the functor of \eqref{FA.dd}, and the {\em
  twisted Hochschild homology object} $CH^F_\idot(\A,\M) \in \D(\E)$
of a small $\dd$-algebra $\A$ with coefficients in an $\A$-bimodule
$\M$ is given by
$$
CH_\idot(\A,\M) = H_\idot(\Delta^o,F_\idot(\M/\A)_\hash),
$$
where $F_\idot(\M/\A)_\hash$ is the functor of \eqref{FAM.dd}. (The
meaning of $H_\idot(\Lambda,-)$ and $H_\idot(\Delta^o,-)$ is the
same as in Subsection~\ref{tw.ho.subs}.)
\item The {\em twisted cyclic homology object} $CC^F_\idot(A_\idot) \in
  \D(\E)$ of a DG algebra $A_\idot$ is given by
$$
CC^F_\idot(\A) = CC^F(\Rr(A_\idot)),
$$
and the {\em twisted Hochschild homology object}
$CH^F_\idot(A_\idot,M_\idot) \in \D(\E)$ of a DG algebra
$A_\idot$ with coefficients in an $A_\idot$-bimodule $M_\idot$ is
given by
$$
CH^F_\idot(\A,\M) = CC^F(\Rr(A_\idot),\Rr(M_\idot)).
$$
\end{enumerate}
\end{defn}

\begin{remark}
If the target category $\E$ is not only exact but also abelian, we
can take the homology of the complexes and define individual twisted
cyclic homology objects $HC^F_i(-) \in \E$, $i \in \Z$, and twisted
Hochschild homology objects $HH_i^F(-,-) \in \E$, $i \in \Z$.
\end{remark}

By definition, $F^{\Delta^o \times \Delta}\A_\hash$ is functorial
with respect to algebra maps, and $F^{\Delta^o \times
  \Delta}(\M/\A)_\hash$ is functorial with respect to maps of pairs
(in both cases, these maps correspond to enriched functors, so that
functoriality follows from \eqref{la.phi}). In other words, for any
map of pairs $f:\langle A_\idot,M_\idot \rangle \to \langle
B_\idot,N_\idot \rangle$ of a DG algebra and a bimodule, we have a
natural map
\begin{equation}\label{F.f.u}
F(f):CH^F_\idot(A_\idot,M_\idot) \to CH^F_\idot(B_\idot,N_\idot),
\end{equation}
and similarly for $\dd$-algebras and bimodules over them.

\begin{lemma}\label{F.exa}
Assume that the trace functor $F:C_\idot(k) \to \E$ is balanced in
the sense of Definition~\ref{adm.fun.def}. Then for any
quasiisomorphism $f:\langle A_\idot,M_\idot \rangle \to \langle
A'_\idot,M'_\idot \rangle$, the corresponding map $F(f)$ is a
quasiisomorphism, and for any total equivalence $f:\langle \A,\M
\rangle \to \langle \A',\M' \rangle$, the corresponding map $F(f)$
is also a quasiisomorphism. In particular, for any $\dd$-algebra
$\A$ and $\A$-bimodule $\M$, we have a natural quasiisomorphism
$$
CH^F_\idot(\A,\M) \cong CH^F(\Ll(\A),\Ll(\M)).
$$
\end{lemma}

\proof{} In the $\dd$-case, note that for a total equivalence
$f:\langle \A,\M \rangle \to \langle \A',M' \rangle$, the
corresponding map
$$
F_\idot(\M/\A)_\hash([n]) \cong F_\idot(\M \otimes \A^{\otimes n-1}) \to
F_\idot(\M'/\A')_\hash([n]) \cong F_\idot(\M' \otimes \A^{'\otimes n-1})
$$
is a quasiisomorphism for any $[n] \in \Delta^o$, so that the
induced map on the homology object $H_\idot(\Delta^o,-)$ is also a
quasiisomorphism. The other statements now follow as in
Lemma~\ref{ba.iso}.
\endproof

We will also need another way to compute the twisted Hochschild
homology objects. Again assume given a bimodule $\M$ over a
$\dd$-algebra $\A$. Then the functor $F^{\Delta^o \times
  \Delta}(\M/\A)_\hash$ of \eqref{Fdd.mu} can be interpreted as an
object
$$
F^{\Delta^o \times \Delta}(\M/\A)_\hash \in (\Delta^o \times \Delta^o
\times \Delta)\E
$$
in the category of functors from the triple product $\Delta^o \times
\Delta^o \times \Delta$ to $\E$. Restricting it with respect to the
diagonal embedding $\delta:\Delta^o \to \Delta^o \times \Delta^o$,
we obtain an object $\delta^*F^{\Delta^o \times \Delta}(\M/\A)_\hash
\in (\Delta^o \times \Delta)\E$. We can then apply the total complex
functor $\Ll$ of \eqref{tot.nd}.

\begin{lemma}\label{F.del}
There exists a natural quasiisomorphism
$$
CH^F_\idot(\A,\M) \cong \Ll(\delta^*F^{\Delta^o \times
  \Delta}(\M/\A)_\hash).
$$
\end{lemma}

\proof{} Since as mentioned in Subsection~\ref{tw.ho.subs}, the
homology $H_\idot(\Delta^o,N)$ of the category $\Delta^o$ with
coefficients in a simplicial object $N$ of some exact category can
be computed by the normalized complex $\Nn_\idot(N)$, the left-hand
side is naturally quasiisomorphic to the total complex of the triple
complex
$$
\Nn^\hdot_{\idot,\idot}(F^{\Delta^o \times \Delta}(\M/\A)_\hash).
$$
To identify it with the right-hand side, it now suffices to apply
the shuffle map \eqref{big.shuffle}.
\endproof

\subsection{DG trace isomorphisms.}

Recall that for any associative unital algebra object $A$ in a
monoidal category $\C$, a right $A$-module $M$, and a left
$A$-module $N$, the {\em bar complex} $\B(M,N/A)$ is a canonical
simplicial object in $\C$ with values given by
$$
\B(M,N/A)([n]) = M \otimes A^{\otimes n-1} \otimes N, \qquad [n] \in
\Delta^o,
$$
and with the structure maps induced by the product in $A$ and by the
$A$-action on $N$ and $M$. If $\C = C_\idot(k)$, this can be applied
to DG modules $M_\idot$, $N_\idot$ over a DG algebra $A_\idot$.
The {\em derived tensor product} $M_\idot \lotimes_{A_\idot}
N_\idot$ is by definition the total complex of the normalization
$$
\Nn_\idot(\B(M_\idot,N_\idot/A_\idot)) \in C_{\idot,\idot}(k).
$$
If we have three DG algebras $A_\idot$, $B_\idot$, $C_\idot$, an
$A_\idot \otimes B_\idot^o$-module $M_\idot$, and a $B_\idot \otimes
C^o_\idot$-module $N_\idot$, then the derived tensor product
$N_\idot \lotimes_{B_\idot} M_\idot$ has a natural structure of an
$A_\idot \otimes C^o_\idot$-module. This construction is functorial
in all the arguments and sends quasiisomorphisms to
quasiisomorphisms.

Assume now given two DG algebras $A_\idot$, $B_\idot$, an $A_\idot
\otimes B^o_\idot$-module $M_\idot$, and a $B_\idot \otimes
A^o_\idot$-module $N_\idot$. What we want to do now is to extend the
arguments of Proposition~\ref{tr.fu.th} and show that for any
balanced trace functor $F:C_\idot(k) \to \E$, the trace functor
structure on $F$ induces a canonical quasiisomorphism
$$
CH^F_\idot(A_\idot, M_\idot \lotimes_{B_\idot} N_\idot) \cong
CH^F_\idot(B_\idot, N_\idot \lotimes_{A_\idot} M_\idot).
$$
We begin with the following observation.

Assume given a $\dd$-algebra $\A$, a right $\A$-module $\M$, and a
left $\A$-module $\N$. Then the bar complex $\B(\M,\N/\A)$ is
naturally an object in $(\Delta^o \times \Delta^o \times
\Delta)k\proj$. Restricting it to the diagonal as in
Lemma~\ref{F.del}, we obtain a object
$$
\overline{\B}(\M,\N/\A) = \delta^*\B(\M,\N/\A) \in (\Delta^o \times
\Delta)k\proj.
$$

\begin{lemma}\label{bar.del}
There exists a natural quasiisomorphism
$$
\Ll(\overline{\B}(\M,\N/\A)) \cong \Ll(\M) \lotimes_{\Ll(\A)} \Ll(\N).
$$
\end{lemma}

\proof{} As in Lemma~\ref{F.del}, the quasiisomorphism is
immediately induced by the shuffle map \eqref{big.shuffle}.
\endproof

We now return to our setting and assume given DG algebras $A_\idot$,
$B_\idot$, an $A_\idot \otimes B^o_\idot$-module $M_\idot$, and a
$B_\idot \otimes A^o_\idot$-module $N\idot$.

\begin{prop}\label{DG.tr}
Under the assumptions above, the trace functor structure on $F$
induces natural isomorphisms
$$
\tau_{M_\idot,N_\idot}:CH^F_\idot(A_\idot, M_\idot
\lotimes_{B_\idot} N_\idot) \cong CH^F_\idot(B_\idot, N_\idot
\lotimes_{A_\idot} M_\idot).
$$
These isomorphisms are functorial and satisfy \eqref{M.N.L}.
\end{prop}

\proof{} Let $\A_1 = \Rr(A_\idot)$, $\A_2 = \Rr(B_\idot)$, $\M_1 =
\Rr(M_\idot)$, $\M_2 = \Rr(N_\idot)$. As in the proof of
Proposition~\ref{tr.fu.th}, consider the two-variable functor
\begin{equation}\label{f.mm}
F^{\Delta^o \times \Delta}(\M_\idot/\A_\idot)_\hash:\Delta^o \times
\Delta^o \to (\Delta^o \times \Delta)\E
\end{equation}
of \eqref{F.multi}. We note that if we restrict it to $[n] \times
\Delta^o$ for some $[n] \in \Delta^o$, then the result coincides on
the nose with
\begin{equation}\label{F.bar}
F^{\Delta^o\times\Delta}(\A_1^{\otimes n-1} \otimes
\B(\M_1,\M_2/\A_2)):\Delta^o \to (\Delta^o \times \Delta)\E,
\end{equation}
where $\B(\M_1,\M_2/\A_2)$ is the bar complex, and $F^{\Delta^o
  \times \Delta}$ is applied pointwise. Let us now interpret
$F^{\Delta^o\times\Delta}(\M_\idot/\A_\idot)_\hash$ as an object
$$
F^{\Delta^o\times\Delta}(\M_\idot/\A_\idot)_\hash \in (\Delta^o
\times \Delta^o \times \Delta^o \times \Delta)\E,
$$
where the first two simplicial directions correspond to the two
simplicial directions in the left-hand side of \eqref{f.mm}, and
denote by $\overline{F(\M_\idot/\A_\idot)}_\hash \in (\Delta^o
\times \Delta)\E$ its restriction to the diagonal $\Delta^o \subset
\Delta^o \times \Delta^o \times \Delta^o$. Then by virtue of
\eqref{F.bar}, Lemma~\ref{bar.del} together with Lemma~\ref{F.del}
provide a canonical quasiisomorphism
$$
CH^F(\Ll(\A_1),\Ll(\M_1) \lotimes_{\Ll(\A_2)} \Ll(\M_2)) \cong
\Ll(\overline{F(\M_\idot/\A_\idot)}_\hash),
$$
and by Lemma~\ref{F.exa}, the left-hand side is canonically
identified with the twisted Hochschild homology object
$CH^F(A_\idot,M_\idot \otimes_{B_\idot} N_\idot)$.

Switching the indices $1$ and $2$, we obtain a canonical
quasiisomorphism
$$
CH^F(B_\idot,N_\idot \lotimes_{A_\idot} M_{\idot}) \cong
\Ll(\overline{F(\M_\idot/\A_\idot)}_\hash);
$$
taken together, these quasiisomorphisms provide the required
quasiisomorphism $\tau_{-,-}$.

To prove \eqref{M.N.L}, one needs to do an obvious three-variable
version of this construction; as in Proposition~\ref{tr.fu.th}, we
leave it to the reader.
\endproof

\subsection{Extended functoriality.}\label{ext.func.subs}

Assume given two admissible $k$-algebras $A$, $B$, and a
multiplicative $k$-linear map $f:A \to B$ which we do {\em not}
assume to be unital. Then $e = f(1) \in B$ is an idempotent element.
Let
$$
\g(f) = eB \in (A,B)\bimod,
$$
with $A$ acting on the left through the map $f$, and let $\g(f)^o =
Be \in (B,A)\bimod$. Then the multiplication map $Be \otimes eB \to
B$ factors through $Be \otimes_A eB$ and induces a $B$-bimodule map
\begin{equation}\label{r.f}
r(f):\g(f) \lotimes_A \g(f)^o \to \g(f) \otimes_A \g(f)^o \to B.
\end{equation}
Moreover, assume given an $A$-bimodule $M$ and an $B$-bimodule
$N$. Then in terms of bimodules $\g(f)$, $\g(f)^o$, extending a map
$f:A \to B$ to a map of pairs $f:\langle A,M \rangle \to \langle B,N
\rangle$ is equivalent to giving an $A$-bimodule map
\begin{equation}\label{l.f}
l(f):M \to \g(f) \lotimes_B N \lotimes_B \g(f)^o \cong \g(f)
\otimes_B N \otimes_B \g(f)^o = eNe \subset N.
\end{equation}
All of the above immediately generalizes to DG algebras and DG
bimodules.

Assume now given a balanced trace functor $\langle
F,\tau_{\idot,\idot} \rangle$ from $k\proj$ to some exact category
$\E$ with infinite filtered colimits. Then for any map of pairs
$$
f:\langle A_\idot,M_\idot \rangle \to \langle B_\idot,N_\idot
\rangle,
$$
Proposition~\ref{DG.tr} provides a quasiisomorphism
$$
CH^F_\idot(A_\idot,\g(f) \lotimes_{B_\idot} N_\idot \lotimes_{B_\idot}
\g(f)^o) \cong CH_\idot^F(B_\idot,\g(f)^o \lotimes_{A_\idot} \g(f)
\lotimes_{B_\idot} N_\idot),
$$
so that \eqref{l.f} induces a map
$$
\begin{CD}
CH_\idot^F(A_\idot,M_\idot) @>{l(f)}>> CH^F(B_\idot,\g(f)^o
  \lotimes_{A_\idot} \g(f) \lotimes_{B_\idot} N_\idot).
\end{CD}
$$
Composing it with the map
$$
\begin{CD}
CH^F_\idot(B_\idot,\g(f)^o \lotimes_{A_\idot} \g(f) \lotimes_{B_\idot}
N_\idot) @>{r(f) \otimes \id}>> CH^F_\idot(B_\idot,N_\idot)
\end{CD}
$$
induced by \eqref{r.f}, we obtain a map
$$
F(f):CH^F_\idot(A_\idot,M_\idot) \longrightarrow
CH^F_\idot(B_\idot,N_\idot).
$$

\begin{prop}\label{funco}
Sending $f$ to $F(f)$ extends $CH^F_\idot(A_\idot,M_\idot)$ to a
functor from the category of pairs $\langle A_\idot,M_\idot \rangle$
of an DG algebra $A_\idot$ over $k$ and an $A_\idot$-bimodule
$M_\idot$, and possibly non-unital maps between them, to the derived
category $\D(\E)$.
\end{prop}

\proof{} We have to check that for any two maps $f:\langle
A_\idot,M_\idot \rangle \to \langle B_\idot,N_\idot\rangle$, $g:
\langle B_\idot,N_\idot \rangle \to \langle C_\idot, L_\idot
\rangle$, we have $F(g \circ f) = F(g) \circ F(f)$, and that if the
map $f$ is unital, then $F(f)$ is the same map as in
\eqref{F.f.u}. Since we obviously have
\begin{equation}\label{compo}
\g(f) \lotimes_{B_\idot} \g(g) \cong \g(g \circ f), \qquad \g(f)^o
\lotimes_{B^o_\idot} \g(g)^o \cong \g(g \circ f)^o,
\end{equation}
the first claim immediately follows from \eqref{M.N.L}. For the
second claim, note that it is tautologically true when $f:A_\idot
\to B_\idot$ is an isomorphism. Therefore we may assume that 
$$
M_\idot \cong \g(f) \lotimes_{B_\idot} N_\idot \lotimes_{B_\idot}
\g(f)^o,
$$
that is, $M_\idot$ is $N_\idot$ considered as an
$A_\idot$-bimodule. The map $F(f)$ of \eqref{F.f.u} is by
definition induced by the map
\begin{equation}\label{f.ddma}
f:F^{\Delta^o \times \Delta}(\M_1/\A_1)_\hash \to F^{\Delta^o \times
  \Delta}(\M_2/\A_2)_\hash,
\end{equation}
where we let $\A_1 = \Rr(A_\idot)$, $\A_2 = \Rr(B_\idot)$, $\M_1 =
\Rr(M_\idot) = \Rr(N_\idot) = \M_2$, as in the proof of
Proposition~\ref{DG.tr}. Moreover, since $\g(f)^o$ considered as a
left $B_\idot$-module is $B_\idot$, $N_\idot \lotimes_{B_\idot}
\g(f)^o$ is naturally homotopy equivalent to $N_\idot$. Therefore,
if we denote $\M_1'=\M_1 = \M_1''=\M_2$, $\A_1'=\A_1$, $\A_1'' =
\A_2$, $\A_2'=\A_2''=\A_2$, $\M_2'=\M_2''=\Rr(\g(f)^o)$, then we can
replace both the source and the target of the map \eqref{f.ddma}
with two-variable functors \eqref{F.multi}, so that $F(f)$ is in
fact induced by the natural map
$$
f:F^{\Delta^o \times \Delta}(\M'_\idot/\A'_\idot)_\hash \to
F^{\Delta^o \times \Delta}(\M''_\idot/\A''_\idot)_\hash.
$$
This map is in turn induced by the natural map of bar complexes
$$
\B(\M'_2,\M'_1/\A_1') \to \B(\M_2',\M'_1/\A''_1),
$$
and this is exactly the adjunction map \eqref{r.f}.
\endproof

\section{DG categories.}\label{dg.cat.sec}

\subsection{Definitions.}

By an {\em admissible DG category} over $k$ we will understand a
small $k$-linear DG category $A_\idot$ such that for any two of its
objects $a$, $a'$, the complex $A_\idot(a,a')$ of maps from $a$ to
$a'$ is a complex of projective $k$-modules. A {\em module} over an
admissible DG category $A_\idot$ is a $k$-linear DG functor $M$ from
$A_\idot$ to $C_\idot(k)$, with the standard DG structure on
$C_\idot(k)$. We denote the category of $A_\idot$-modules by
$A_\idot\amod$. Inverting pointwise quasiisomorphisms in $A_\idot$
gives the derived category $\D(A_\idot)$.

For an admissible DG category $A_\idot$, we denote by $A^o_\idot$
the opposite category. For two admissible DG categories $A_\idot$,
$B_\idot$ over $k$, we denote by $A_\idot \otimes B_\idot$ their
tensor product -- that is, the DG category of pairs $\langle a,b
\rangle$ of an object $a \in A_\idot$ and an object $b \in B_\idot$,
with maps given by
$$
(A_\idot \otimes B_\idot)(\langle a,b \rangle, \langle a',b'
\rangle) = A_\idot(a,a') \otimes B_\idot(a,b').
$$
We denote $A_\idot\bimod = (A^o_\idot \otimes A_\idot)\amod$;
objects of this category are {\em $A_\idot$-bimodules}.

Assume given an admissible DG category $A_\idot$ over $k$. Then for
any finite set $S$ of objects in $A_\idot$, we define a DG algebra
$A^S_\idot$ over $k$ by
$$
A^S_\idot = \bigoplus_{a,a' \in S}A_\idot(a,a'),
$$
with the obvious multiplication and the unity element. The DG
algebra $A^S_\idot$ is admissible. For any $V_\idot \in
A_\idot\amod$, we let
$$
V_\idot^S = \bigoplus_{a \in S}V_\idot(a),
$$
and for any $M_\idot \in A_\idot\bimod$, we let
$$
M_\idot^S = \bigoplus_{a,a' \in S}M_\idot(a,a').
$$
We obviously have $V^S_\idot \in A^S_\idot\amod$, $M^S_\idot \in
A^S_\idot\bimod$. Moreover, for every embedding $S_1 \to S_2$ of finite
sets of objects in $A_\idot$, and for any $M_\idot \in
A_\idot\bimod$, we have a natural map of pairs
\begin{equation}\label{incl}
\langle A^{S_1}_\idot,M^{S_1}_\idot \rangle \to \langle
A^{S_2}_\idot,M^{S_2}_\idot \rangle.
\end{equation}
Assume now given a balanced trace functor $\langle
F,\tau_{\idot,\idot}\rangle$ from $k\proj$ to a Karou\-bi-closed exact
category $\E$ with arbitrary sums. Then by Proposition~\ref{funco},
the maps \eqref{incl} induce natural maps
$$
CH^F_\idot(A^{S_1}_\idot,M^{S_1}_\idot) \to
CH^F_\idot(A^{S_2}_\idot,M^{S_2}_\idot).
$$

\begin{defn}
For any admissible DG category $A_\idot$ over $k$ and any $M_\idot
\in A_\idot\bimod$, we have
$$
CH_\idot^F(A_\idot,M_\idot) =
\lim_{\overset{S}{\to}}CH^F_\idot(A^S_\idot,M^S_\idot),
$$
where the limit is taken over all finite sets of objects in
$A_\idot$.
\end{defn}

We note that since the limit is filtered, it descends to the derived
category level with no ambiguity.

Moreover, assume given three admissible DG categories $A_\idot$,
$B_\idot$, $C_\idot$ over $k$, an $A_\idot^o \otimes B_\idot$-module
$M_\idot$, and a $B^o_\idot \otimes C_\idot$-module $N_\idot$. Then
for any objects $a$ in $A_\idot$ and $c$ in $C_\idot$, the filtered
direct limit
$$
M_\idot \lotimes_{B_\idot} N_\idot(a,c) =
\lim_{\overset{S}{\to}}M^{\{a\} \times S}_\idot \lotimes_{B^S_\idot}
N^{S \times \{c\}}_\idot
$$
over all finite sets of objects in $B_\idot$ is well-defined, and
gives a natural tensor product module
$$
M_\idot \lotimes_{B_\idot} N_\idot \in (A_\idot^o \otimes
C_\idot)\amod.
$$
Then if $C_\idot = A_\idot$, Proposition~\ref{DG.tr} provides
natural quasiisomorphisms
\begin{equation}\label{cat.tr}
\tau_{M_\idot,N_\idot}:CH^F_\idot(A_\idot, M_\idot \lotimes_{B_\idot}
N_\idot) \cong CH^F_\idot(B_\idot, M_\idot \lotimes_{A_\idot}
N_\idot)
\end{equation}
satisfying \eqref{M.N.L}.

Assume now given a DG functor $f:A_\idot \to B_\idot$ between
admissible DG categories over $k$. Then $f$ induces a pullback
functor $f^*:B_\idot\amod \to A_\idot\amod$ descending to a functor
$f^*:\D(A_\idot) \to \D(B_\idot)$. In terms of tensor products,
$f^*$ is given by
$$
f^*N_\idot \cong \g(f) \lotimes_{B_\idot} N_\idot,
$$
where $\g(f) \in (A_\idot^o \otimes B_\idot)\amod$ is the graph of
the functor $f$ given by
\begin{equation}\label{g.f.eq}
\g(f)(a,b) = \B_\idot(f(a),b).
\end{equation}
The functor $f$ also defines a functor $f^o:A_\idot^o \to
B_\idot^o$, and we denote $\g(f)^o = \g(f^o)$. Note that this is
consistent with our earlier notation -- for any finite set $S$ of
objects in $A_\idot$, we have
$$
\g(f)^S \cong \g(f^S), \g(f)^{oS} \cong \g(f^S)^o
$$
where $f^S:A_\idot^S \to B_\idot^{f(S)}$ is the natural map induced
by the functor $f$. Setting
$$
f_!(M_\idot) = \g(f)^o \otimes_{A_\idot} M_\idot
$$
defines a functor $f_!:\D(A_\idot) \to \D(B_\idot)$ left-adjoint to
the pullback functor $f^*$. The natural map
$$
r(f):\g(f)^o \lotimes_{A_\idot} \g(f) \to B_\idot \cong \g(\id)
$$
of \eqref{l.f} induces the adjunction. For any two composable
functors $f:A_\idot \to B_\idot$, $g:B_\idot \to C_\idot$, we have
$(g \circ f)^* \cong f^* \circ g^*$ and
$$
\g(f) \lotimes_{B_\idot} \g(g) \cong \g(g \circ f),
$$
as in \eqref{compo}.

The functor $f$ also extends to a functor $f:A^o_\idot \otimes
A_\idot \to B^o_\idot \otimes B_\idot$, so that we obtain a pullback
functor $f^*:A_\idot\bimod \to B_\idot\bimod$. Given an
$A_\idot$-bimodule $M_\idot$ and an $B_\idot$-bimodule $N_\idot$, we
will say that $f$ is extended to a {\em map of pairs} $\langle
A_\idot,M_\idot \rangle \to \langle B_\idot,N_\idot \rangle$ if we
are given an $A_\idot$-bimodule map
$$
l(f):M_\idot \to f^*N_\idot \cong \g(f) \lotimes_{B_\idot} N_\idot
\lotimes_{B_\idot} \g(f)^o.
$$
Then as in Subsection~\ref{ext.func.subs}, the maps $l(f)$, $r(f)$
together with the trace isomorphisms \eqref{cat.tr} give a natural
map
$$
CH_\idot^F(A_\idot,M_\idot) \to CH_\idot^F(B_\idot,N_\idot)
$$
for any map $f:\langle A_\idot,M_\idot \rangle \to \langle
B_\idot,N_\idot \rangle$, and these maps are associative in the
obvious sense (for the proof, one takes the corresponding staements
for DG algebras and passes to the appropriate filtered direct
limit). In particular, we always have a natural map
\begin{equation}\label{CH.f}
CH_\idot^F(A_\idot,f^*N_\idot) \to CH_\idot^F(B_\idot,N_\idot);
\end{equation}
the corresponding map $l(f)$ is the identity map. Another canonical
extension is a map $f:\langle A_\idot,A_\idot \rangle \to \langle
B_\idot,B_\idot \rangle$; this gives a natural map
\begin{equation}\label{CH.f.2}
CH_\idot^F(A_\idot) \to CH_\idot^F(B_\idot).
\end{equation}

\begin{prop}[``Density'']\label{dens}
Assume given a functor $f:A_\idot \to B_\idot$ between admissible DG
categories such that $f^*:\D(B_\idot) \to \D(A_\idot)$ is a fully
faithful embedding, and a balanced trace functor $F:k\proj \to \E$
to some Karoubi-closed exact category with arbitrary sums. Then for
any $B_\idot$-bimodule $N_\idot$, the map \eqref{CH.f} is a
quasiisomorphism.
\end{prop}

\proof{} Since $f^*$ is fully faithful, we have $f_! \circ f^* \cong
\id$, so that
$$
\g(f)^o \lotimes_{A_\idot} \g(f) \cong B_\idot,
$$
and the map $r(f)$ is an isomorphism. Since the map $l(f)$ is the
identity map, we are done.
\endproof

\begin{corr}[``Morita-invariance'']
In the assumption of Proposi\-ti\-on~\ref{dens}, assume in addition
that $f^*$ is an equivalence of categories. Then the natural map
\eqref{CH.f.2} is an isomorphism.
\end{corr}

\proof{} The condition ensures that $f^* \circ f_! \cong \id$ as
well, so that the map $l(f)$ is an isomorphism.
\endproof

In particular, for any admissible DG category $A_\idot$ over $k$,
its pretriangulated hull $\wt{A}_\idot$ is also an admissible DG
category over $k$, and we have a canonical identification
\begin{equation}\label{hull}
CH_\idot^F(A_\idot) \cong CH_\idot^F(\wt{A}_\idot).
\end{equation}

\subsection{Localization.}\label{loc.subs}

Assume now given a sequence
\begin{equation}\label{seq}
\begin{CD}
A_\idot @>{f}>> B_\idot @>{g}>> C_\idot
\end{CD}
\end{equation}
of admissible DG categories over $k$ and DG functors between them.

\begin{defn}
A sequence \eqref{seq} is {\em exact} if $f^* \circ g^* = 0$, $g^*$
is a fully faithful embedding, and $f^*$ induces an equivalence
$$
\D(B_\idot)/g^*(\D(C_\idot)) \cong \D(A_\idot).
$$
\end{defn}

As we have mentioned in the introduction, Keller's Theorem \cite{kel}
asserts that for an exact sequence \eqref{seq}, we have a long exact
sequence
$$
\begin{CD}
HH_\idot(A_\idot) @>>> HH_\idot(B_\idot) @>>> HH_\idot(C_\idot) @>>>
\end{CD}
$$
of Hochschild homology groups. This is very easy to see in our
language. Indeed, the assumption of exactness implies that we have
$f^* \circ f_! \cong \id$ and $g_! \circ g^* \cong \id$, so that
\begin{equation}\label{tr.exa}
\begin{aligned}
HH_\idot(A_\idot) &\cong HH_\idot(B_\idot,\g(f)^o \lotimes_{A_\idot}
\g(f)),\\
HH_\idot(C_\idot) &\cong HH_\idot(B_\idot,\g(g) \lotimes_{C_\idot}
\g(g)^o),
\end{aligned}
\end{equation}
where $\g(f)$, $\g(g)$ are as in \eqref{g.f.eq}. Moreover, again by
exactness, we have a distinguished triangle
\begin{equation}\label{tria.1}
\begin{CD}
\g(f)^o \lotimes_{A_\idot} \g(f) @>>> B_\idot @>>> \g(g)
\lotimes_{C_\idot} \g(g)^o @>>>
\end{CD}
\end{equation}
of $B_\idot$-bimodules, and since the Hochschild homology functor
$HH_\idot(B_\idot,-)$ is exact, this induces the desired long exact
sequence. We note that the identifications \eqref{tr.exa} hold for
twisted Hochschild homology groups just as well. This motivates the
following definition.

\begin{defn}\label{loca}
A balanced trace functor $F:k\proj \to \E$ is {\em localizing} for
any exact sequence \eqref{seq} of DG categories, the triangle
$$
\begin{CD}
CH_\idot^F(B_\idot,\g(f)^o \lotimes_{A_\idot} \g(f)) @>>>
CH_\idot^F(B_\idot) @>>>\\
@>>> CH_\idot^F(B_\idot,\g(g) \lotimes_{C_\idot} \g(g)^o) @>>>
\end{CD}
$$
in $\D(\E)$ induced by \eqref{tria.1} is a distinguished triangle.
\end{defn}

Thus for a localizing balanced trace functor, we have a distiguished
triangle
$$
\begin{CD}
CH^F_\idot(A_\idot) @>>> CH^F_\idot(B_\idot) @>>> CH^F_\idot(C_\idot) @>>>
\end{CD}
$$
for any exact sequence \eqref{seq}.

\begin{lemma}
An extension and a filtered direct limit of localizing balanced
trace functors is localizing.
\end{lemma}

\proof{} Clear. \endproof

We note that there is no hope that an arbitrary balanced trace
functor is localizing: for example, if we take a balanced functor
$F:k\proj \to \E$ and equip it with the trivial trace functor
structure of Example~\ref{triv.exa}, then we have
$$
CH_\idot^F(A_\idot) \cong F(CH_\idot(A_\idot)),
$$
and $F$ has no hope to be localizing unless it is additive. We
finish the paper with proving that there is at least one non-trivial
example: the cyclic power trace functor of Example~\ref{cycl.exa} is
localizing in the sense of Definition~\ref{loca}.

\begin{lemma}\label{cycl.fil}
For any $l \geq 1$, the cyclic power functor $\Cycl^l:k\proj \to
k\proj$ of Example~\ref{cycl.exa} is filterable in the sense of
Definition~\ref{bounded.defn}.
\end{lemma}

\proof{} Equip the $l$-th tensor power $M^{\otimes_k l}$ of some
$M_\idot \in W_\idot k\proj$ with a grading in the usual way,
namely,
\begin{equation}\label{gra.cy}
\left(M^{\otimes_k l}\right)_i = \bigoplus_{\rk(f)=i} M_{f(1)}
\otimes \dots \otimes M_{f(l)},
\end{equation}
where the sum is over all maps $f$ from the set $1,\dots,l$ to the
set of non-negative integers such that
$$
\rk(f) = \sum_{1 \leq j \leq l}f(j)
$$
is equal to $i$. The action of the cyclic group $\Z/l\Z$ preserves
this grading, so that the grading descends to the quotient
$\Cycl^l(M)$. The corresponding filtration is obviously functorial
with respect to filtered maps $M_\idot \to M'_\idot$ and satisfies
the conditions of Definition~\ref{bounded.defn}.
\endproof

Thus by Proposition~\ref{fil=>bal}, $\Cycl^l$ is also balanced in
the sense of Definition~\ref{adm.fun.def}, so that
$CH^{\Cycl^l}_\idot(-,-)$ is well-behaved. To simplify notation, let
$$
CH^{(l)}_\idot = CH^{\Cycl^l}_\idot
$$
for any $l \geq 1$. For any integer $n \geq 0$, let $H(n)$ be the
set of all maps $f:\{1,\dots,l\} \to \{0,\dots,n\}$. The cyclic
group $\Z/l\Z$ acts on $H(n)$ via its action on the set
$\{1,\dots,l\}$. As in the proof of Lemma~\ref{cycl.fil}, let
$$
\rk(f) = \sum_{1 \leq i \leq l}f(i)
$$
for any $f \in H(n)$, and let $p(f)$ be the order of the stabilizer
of $f$ in $\Z/l\Z$. We then have $l = p(f)q(f)$ for some integer
$q(f)$, and we have $f(i+q(f)) = f(i)$ for any $i$, $1 \leq i \leq
l-q(f)$. For any $i$, $0 \leq i \leq l$, let $H(n,i) \subset H(n)$
be the subset of all maps $f$ with $\rk(f)=i$, and let $Q(n,i) =
H(n,i)/(\Z/l\Z)$. Fix a splitting $Q(n,i) \to H(n,i)$ of the
quotient map $H(n,i) \to Q(n,i)$.

\begin{lemma}\label{cycl.A.dec}
Assume an admissible $k$-algebra $A$ and projective $A$-bimodules
$M_0,\dots,M_n$. We then have a functorial decomposition
$$
\Cycl^l_A(M_0 \oplus \dots \oplus M_n) = \bigoplus_{0 \leq i \leq
  l}\Cycl^{l,i}_A(M_0,\dots,M_n),
$$
where
\begin{equation}\label{deco.0}
\Cycl^{l,i}_A(M_0,\dots,M_n) =
\bigoplus_{f \in Q(n,i)}\Cycl^{p(f)}_A(M_{f(1)} \otimes_A \dots
\otimes_A M_{f(q(f))}).
\end{equation}
\end{lemma}

\begin{remark}
Formally, the right-hand side of the decomposition \eqref{deco.0}
depends on the choice of a splitting $Q(n,i) \to H(n,i)$; however,
objects corresponding to different choices of splittings are
canonically identified by the trace isomorphisms
$\tau_{\idot,\idot}$.
\end{remark}

\proof{} Use the explicit description of the functor $\Cycl^l_A$
given in Lemma~\ref{cycl.A}. Then \eqref{gra.cy} induces the desired
decomposition, and to identify the terms, it remains to take the
quotient with respect to the cyclic group.
\endproof

\begin{prop}
For any $l \geq 1$, the balanced trace functor $\Cycl^l$ of
Example~\ref{cycl.exa} is localizing in the sense of
Definition~\ref{loca}.
\end{prop}

\proof{} Applying Lemma~\ref{cycl.A.dec} pointwise to
simplicial-cosimplicial objects, we see that for any
$B_\idot$-bimodules $M_0,\dots,M_n$, we have a natural decomposition
$$
CH^{(l)}_\idot(B_\idot,M_0 \oplus \dots \oplus M_n) =
\bigoplus_{0 \leq i \leq l}CH^{(l),i}_\idot(B_\idot,M_0,\dots,M_n)
$$
with
\begin{multline}\label{deco}
CH^{(l),i}_\idot(B_\idot,M_0,\dots,M_n) = \\
=\bigoplus_{f \in Q(n,i)}CH^{p(f)}_\idot(B_\idot,M_{f(1)}
\lotimes_{B_\idot} \dots \lotimes_{B_\idot} M_{f(q(f))}).
\end{multline}
Take $n=1$, $M_0 = \g(f)^o \lotimes_{A_\idot} \g(f)$, $M_1 = \g(g)
\lotimes_{C_\idot} \g(g)^o$. Then since the sequence \eqref{seq} is
exact, we have
$$
M_0 \lotimes_{B_\idot} M_1 = M_1 \lotimes_{B_\idot} M_0 = 0.
$$
Therefore the only non-trivial terms in \eqref{deco} are those with
$q(f) = 1$, that is, $f:\{1,\dots,l\} \to \{0,1\}$ constant. Thus 
$$
CH^{(l),i}(B_\idot,M_0, M_1) = 0
$$
unless $i=0$ or $i=l$, and in these cases, we have
$$
\begin{aligned}
CH^{(l),0}_\idot(B_\idot,M_0, M_1) &\cong
CH^{(l)}_\idot(B_\idot,M_0),\\
CH^{(l),l}_\idot(B_\idot,M_0, M_1) &\cong
CH^{(l)}_\idot(B_\idot,M_1).
\end{aligned}
$$
Choosing appripriate DG representatives for $M_0$, $M_1$, $M \cong
B_\idot$, we can assume that the triangle \eqref{tria.1} is actually
a short exact sequence of $B_\idot$-bimodules, or equivalently, a
two-term increasing filtration on a $B_\idot$-bimodule $M$
quasiisomorphic to $B_\idot$. By Lemma~\ref{cycl.fil}, this
filtration induces a generalized filtration $W_\idot$ on
$CH^{(l)}(B_\idot,M)$ with
$$
\gr^i_WCH_\idot^{(l)}(B_\idot,M) \cong
CH^{(l),i}_\idot(B_\idot,M_0,M_1).
$$
We conclude that $\gr^i_WCH^{(l)}_\idot(B_\idot,M)=0$ unless $i=0,l$,
and the triangle
$$
\begin{CD}
CH^{(l)}_\idot(B_\idot,M_0) @>>> CH^{(l)}(B_\idot,M) @>>>
CH^{(l)}_\idot(B_\idot,M_1) @>>>
\end{CD}
$$
induced by the generalized filtration $W_\idot$ is
distinguished. This is what we had to prove.
\endproof

\bigskip

\noindent
{\sc
Steklov Math Institute, Algebraic Geometry section\\
\mbox{}\hspace{30mm}and\\
Laboratory of Algebraic Geometry, NRU HSE}

\bigskip

\noindent
{\em E-mail address\/}: {\tt kaledin@mi.ras.ru}

\end{document}